\documentclass{article}

\input{zstyle/preamble}

\pagecolor{white}

\begin{document}

\maketitle

\begin{abstract} \noindent
Using profinite Galois descent, we compute the Brauer group of the $K(1)$-local category relative to Morava E-theory. At odd primes this group is generated by a cyclic algebra formed using any primitive $(p-1)$st root of unity, but at the prime two is a group of order $32$ with nontrivial extensions; we give explicit descriptions of the generators, and consider their images in the Brauer group of $KO$. Along the way, we compute the relative Brauer group of completed $KO$, using the \'etale locally trivial Brauer group of Antieau, Meier and Stojanoska.
\end{abstract}

\tableofcontents


\section{Introduction}
\label{sec:intro}

The classification of central simple algebras over a field is a classical question in number theory, and by the Wedderburn theorem any such algebra is a matrix algebra over some division algebra, determined up to isomorphism. If one identifies those algebras which arise as matrix algebras over the same division ring, then tensor product defines a group structure on the resulting set of equivalence classes. This relation is \emph{Morita equivalence}, and the resulting group is the \emph{Brauer group} $\Br(K)$. One formulation of class field theory is the determination of $\Br(K)$ in the case that $K$ is a number field.

One consequence of Wedderburn's theorem is that every central simple algebra is split by some extension $L/K$; in fact, one can take $L$ to be Galois. This opens the door to cohomogical descriptions of the Brauer group: by Galois descent, one obtains the first isomorphism in
    \[ \Br(K) \cong H^1(K, \PGL_\infty) \cong H^2(K, \G_m), \label{intro:Br = H2 Gm} \numberthis \]
where $\PGL_\infty$ denotes the Galois module $\varinjlim \PGL_k(K^{\mathrm{s}})$. The second isomorphism is \cite[Proposition~X.9]{serre_localfields}, and follows from Hilbert 90. This presents $\Br(K)$ as an \'etale cohomology group, and allows the use of cohomological techniques in its determination. Conversely, it gives a concrete interpretation of $2$-cocycles, analogous to the relation between $1$-cocycles and Picard elements.

Globalising this picture was a deep problem in algebraic geometric, initiated by the work of Azumaya, Auslander and Goldman, and the Grothendieck school. There are two ways to proceed: on one hand, one can generalise the notion of central simple algebra, obtaining that of an \emph{Azumaya algebra} on $X$, and the group $\Br(X)$ classifying such algebras (again, up to Morita equivalence). On the other, one can use the right-hand side of \cref{intro:Br = H2 Gm}, which globalises in an evident way; this yields the more computable group $\Br'(X) \coloneqq H^2(X, \G_m)_{\mathrm{tors}}$. A famous problem posed by Grothendieck asks if these groups agree in general. While there is always an injective map $\Br(X) \hookrightarrow \Br'(X)$, surjectivity is known to fail for arbitrary schemes (e.g. \cite[Corollary~3.11]{edidin_brauer_2001}). The fundamental insight of To\"en was that it is key to pass to \emph{derived} Azumaya algebras, in which case one \emph{can} represent any class in $H^2(X, \G_m)$ as a derived Azumaya algebra, including nontorsion classes \cite{toen}. To\"en's work shows that even the Brauer groups of classical rings are most naturally studied in the context of derived or homotopical algebraic geometry. This initiated a study of Brauer groups through the techniques of higher algebra, and in recent years they have become objects of intense study in homotopy theory.

In this context, the basic definitions appear in \cite{baker-richter-szymik}: as in the classical case, the Brauer group of a ring spectrum $R$ classifies Morita equivalence classes of Azumaya algebras, defined by entirely analogous axioms. This gives a useful class of noncommutative algebras over an \Einfty-ring: for example, Hopkins and Lurie \cite{hl} computed Brauer groups as part of a program to classify $E_h$-algebra structures on Morava K-theory, where $E_h$ denotes Lubin-Tate theory. Categorifying, the Azumaya condition on an $R$-algebra $A$ is equivalent to the requirement that the \inftycat{} $\Mod_A$ be $R$-linearly invertible, and so the Brauer group is also the Picard group of the \Einfty-algebra $\Mod_R \in \Pr^L$; in this guise the Brauer space and its higher-categorical analogues are the targets for (invertible) factorisation homology, as an extended quantum field theory \cite{lurie_cobordism,scheimbauer_thesis,ayala-francis,haugseng_morita,freed-hopkins_reflectionpositivity}. Computations of Brauer groups of ring spectra of particular interest appear in \cite{antieau-gepner,gl,antieau-meier-stojanoska}.

In this document we study Brauer groups in the monochromatic setting. Recall that Hopkins and Lurie focus on the $K(h)$-local Brauer group of $E_h$. Our objective is to extend this to the Brauer group of the entire $K(h)$-local category: explicitly, this group classifies of $K(h)$-local Azumaya algebras up to Morita equivalence. The computations of \opcit\ show that the groups $\Br(E_h)$ are highly nontrivial, unlike the analogous Picard groups. We therefore complement this by computing the \emph{relative} Brauer group, which classifies those $K(h)$-local Azumaya algebras which become trivial after basechange to $E_h$. In analogy with the established notation for Picard groups (e.g. \cite[p.5]{ghmr_pic}), this group will be denoted $\Br_h^0$. It fits in an exact sequence
    \[ 1 \to \Br_h^0 \to \Br_h \to \Br(E_h)^{\G_h}, \]
where $\G_h$ denotes the Morava stabiliser group, and this gives access (at least in theory) to the group $\Br_h \coloneqq \Br(\mathcal Sp_{K(h)})$. The group $\Br_h^0$ also has a concrete interpretation in terms of chromatic homotopy theory: it classifies twists of the $\G_h$-action on the \inftycat{} $\Mod_{E_h}(\mathcal Sp_{K(h)})$. For the standard action (by basechange along the Goerss-Hopkins-Miller action on $E_h$), one has $\mathcal Sp_{K(h)} \simeq \Mod_{E_h}(\mathcal Sp_{K(h)})^{h\G_h}$ (see \cite[Proposition~10.10]{mathew_galois} and \cite[Theorem~A.II]{profinitedescent} for two formulations). Taking fixed points for a twisted action therefore gives a twisted version of the $K(h)$-local category.

Our main theorems give the computation of the relative Brauer group at height one:

\begin{introthm}[\cref{lem:Brauer p = 2}] \label{intro:Br1 even}
At the prime two, 
    \[ \Br_1^0 \cong \Z/8 \{Q_1\} \oplus \Z/4 \{Q_2\}. \]
        
\begin{enumerate}
    \item The $\Z/4$-factor is mapped injectively to $\Br(KO_2)$, and $Q_2^{\otimes 2} \otimes KO_2$ is the image of the generator under
        \[ \Z/2 \cong \Br(KO \mid KU) \to \Br(KO_2 \mid KU_2). \]
    The $\Z/8$ factor is the relative Brauer group $\Br(\mathcal Sp_{K(1)} \mid KO_2)$.
    
    \item $Q_4 \coloneqq Q_1^{\otimes 4}$ is the cyclic algebra formed using the $C_2$-Galois extension $\mathbb S_{K(1)} \to KU_2^{h (1 + 4\Z[2])}$ and the \emph{strict} unit
        \[ 1 + \varepsilon \in \pi_0 \GL_1(\mathbb S_{K(1)}) = (\Z[2][\varepsilon]/(2 \varepsilon, \varepsilon^2))^\times. \]
\end{enumerate}
\end{introthm}

We have indexed generators on the filtration in which they are detected in the descent spectral sequence, which we recall later in the introduction. For the final part, note \cite[§4]{baker-richter-szymik} that cyclic algebras are defined using \emph{strict} units: that is, maps of spectra
    \[ u : \Z \to \gl_1(E). \]
We give an alternative construction of cyclic algebras from strict units in \cref{sec:cyclic}, and using this we show that they are detected in the HFPSS by a \emph{symbol} in the sense of \cite[Chapter~XIV]{serre_localfields}. This allows us to deduce when they give rise to nontrivial Brauer classes.

Any strict unit has an underlying unit, and we abusively denote these by the same symbol. The unit $1 + \varepsilon$ was shown to be strict in \cite{carmeli-yuan}, which is what is what gives rise to the claimed representative for the class $Q_4$ at the prime two. Likewise, at odd primes the roots of unity are strict, which leads to our second main computation:

\begin{introthm}[\cref{lem:Brauer p odd}]
\label{intro:Br1 odd}
At odd primes,
    \[ \Br_1^0 \cong \Z/(p-1). \] 
A generator is given by the cyclic algebra $(KU_p^{h(1 + p\Z[p])}, \chi, \omega)$, where $\chi: \mu_{p-1} \cong \Z/(p-1)$ is a character and $\omega \in (\pi_0 \mathbb S_{K(1)})^\times \cong \Z[p]^\times$ is a primitive $(p-1)$st root of unity.
\end{introthm}

We now give an outline of the computation. Since Grothendieck, the main approach to computing Brauer groups has been \emph{\'etale} or \emph{Galois descent}, and this is the case for us too. Namely, recall that at any height $h$, Morava E-theory defines a $K(h)$-local Galois extension
    \[ \mathbb S_{K(h)} \to E_h, \]
with profinite Galois group $\G_h$. In \cite{profinitedescent}, we used condensed mathematics to prove a Galois descent statement of the form
    \[ \mathcal Sp_{K(h)} \simeq \Mod_{E_h}(\mathcal Sp_{K(h)})^{h\G_h}, \]
and deduced from this a homotopy fixed point spectral sequence for Picard and Brauer groups, extending the Galois descent results of Mathew and Stojanoska \cite{ms} and Gepner and Lawson \cite{gl}. For our purposes, the main computational upshot of that paper is:

\begin{thm}
There is a descent spectral sequence
        \[ E_2^{s,t} = H^s(\G_h, \pi_t \pic(E_h)) \implies \pi_{t-s} \pic(\mathcal Sp_{K(h)}). \]
whose $(-1)$-stem gives an upper bound on $\Br_h^0$. In a large range, there is an explicit comparison of differentials with the $K(h)$-local $E_h$-Adams spectral sequence.
\end{thm}

A more precise form of the theorem is recalled in \cref{sec:sseq}. In the present paper, we determine this spectral sequence completely at height one, at least for $t-s \geq -1$. Unsurprisingly, the computation looks very different in the cases $p = 2$ and $p > 2$, and the former represents the majority of our work. To prove a lower bound on $\Br_1^0$, we also prove a realisation result in the spirit of To\"en's theorem. Namely, we show in a very general context that all classes on the $E_\infty$-page may be represented by Azumaya algebras:

\begin{introthm}[\cref{lem:Br = Br'}] \label{intro:Br = Br'}
Let $\mathcal C$ be a nice symmetric monoidal \inftycat{} and $A \in \CAlg(\mathcal C)$ a faithful and dualisable Galois extension of the unit. Then the map sending an Azumaya algebra to its module \inftycat{},
    \[ \Br(\mathcal C \mid A) \to \Br'(\mathcal C \mid A) \coloneqq \{\mathcal D : \Mod_A (\mathcal D) \simeq \Mod_A(\mathcal C) \} \subset \Pic(\Mod_{\mathcal C}({\Pr}^L)), \]
is an isomorphism.
\end{introthm}

See also \cite[§6.3]{antieau-gepner} and \cite[§6.4]{gl} for related results. We refer the reader to \cref{sec:generators} for the exact conditions in the theorem; most pertinently, the chromatic localisations of spectra give examples of such \inftycats{}. This accounts immediately for most of the classes in $\Br_1^0$, by sparsity in the descent spectral sequence. At the prime two there is one final computation necessary, which is the relative Brauer group of $KO_2$. Recall that the group $\Br(KO \mid KU)$ was computed by Gepner and Lawson \cite[]{gl}. In \cref{sec:K-theory}, we prove:

\begin{introthm}[\cref{lem:Brauer of KO2}] \label{intro:Br(KO2 mid KU2)}
The relative Brauer group of $KO_2$ is
    \[ \Br(KO_2 \mid KU_2) \cong \Z/4, \]
and the basechange map from $\Br(KO \mid KU) \cong \Z/2$ is injective.
\end{introthm}

The generator, which we denote $P_2$, may be thought of as a cyclic algebra for the unit $1 + 4\zeta$, where $\zeta$ is a topological generator of $\Z[2]$ (see \cref{sec:cyclic}). We remark that this unit is \emph{not} strict, and so the cyclic algebra cannot be constructed by hand; instead, we invoke \cref{intro:Br = Br'} to prove that the possible obstructions vanish. We then show that $P_2$ survives the descent spectral sequence for the extension $\mathbb S_{K(1)} \to KO_2$, and therefore gives rise to the class $Q_2 \in \Br_1^0$.

As an aside, we observe that \cref{intro:Br(KO2 mid KU2)} implies the following result, which may be of independent interest:

\begin{introthm} \label{intro:splitting KUp}
There is no $C_2$-equivariant splitting
    \[ \gl_1 KU_2 \simeq \tau_{\leq 3} \gl_1 KU_2 \oplus \tau_{\geq 4} \gl_1 KU_2. \]
\end{introthm}

We do not know if such an equivariant splitting exists for $\gl_1 KU$. It seems conceivable that one would arise from a discrete model for K-theory, and this is an interesting problem.

\subsection{Outline of the paper}
In \cref{sec:sseq}, we compute the $E_2$-page of the descent spectral sequences and most differentials, giving an upper bound on the relative Brauer groups. The objective of the rest of the document is to show this bound is achieved. In \cref{sec:generators} we prove that the property of ``admitting a compact generator'' satisfies $\mathcal Sp_{K(h)}$-linear Galois descent, which allows us to lift certain $E_\infty$-classes to Azumaya algebras. In \cref{sec:cyclic} we give a construction of cyclic algebras which makes clear where they are detected in the descent spectral sequence; this allows us to assert that the cyclic algebras we form are distinct, and completes the odd-prime computation. We also give a construction of Brauer classes from $1$-cocycles, which we use at the prime two. In \cref{sec:azumaya} we compute the relative Brauer group $\Br(KO_2 \mid KU_2)$, and use this to complete the computation of $\Br_1^0$ at the prime two.

\subsection{Acknowledgements}
I'm extremely grateful to my PhD supervisors, Behrang Noohi and Lennart Meier, for their support throughout the past few years; I'd also like to thank them for useful comments on an earlier draft. I also had helful conversations on this project with Zachary Halladay, Kaif Hilman, Shai Keidar and Sven van Nigtevecht. I'd especially like to thank Maxime Ramzi for suggesting the strategy of \cref{lem:widehat u is E1}, and for comments on a previous version. This work forms part of my thesis, supported by EPSRC [grant EP/R513106/1].

\subsection{Notation and conventions}
\begin{itemize}
\item We will freely use the language of \inftycats{} (modeled as quasi-categories) as pioneered by Joyal and Lurie \cite{htt, ha, sag}. In particular, all (co)limits are $\infty$-categorical. Most commonly, we will be in the context of a presentably symmetric monoidal stable \inftycat{}, and we use the term \emph{stable homotopy theory} to mean such an object. All our computations take place internally to the $K(h)$-local category, and so the symbol $\otimes$ will  generally denote the $K(h)$-local smash product.

\item We only consider spectra with group actions, and not any more sophisticated equivariant notion. When $G$ is a profinite group, we will write $H^*(G,M)$ for \emph{continuous} group cohomology.

\item We follow the conventions of \cite{profinitedescent}. In particular, we direct the reader to \cite[§2.2]{profinitedescent} for details about the pro\'etale site and the sheaves $\mathcal E$ and $\pic(\mathcal E)$. Given a descendable $G$-Galois extension $\bm 1 \to A$ in a stable homotopy theory (possibly profinite), we will implicitly use the associated (hypercomplete) sheaf $\mathcal A \in \Sh(\proet G, \mathcal C)$, writing $A^{hG} \coloneqq \Gamma \mathcal A$. Using \cite[§3.1]{profinitedescent}, we will also form the sheaf $\Mod_{\mathcal A}(\mathcal C) \in \Sh(\proet G, \Pr^{L, \mathrm{smon}})$, and hence the Picard sheaf $\pic(\mathcal A) \in \Sh(\proet G, \mathcal Sp_{\geq 0})$. In this case the descent spectral sequence reads
    \[ H^s(\proet G, \pi_t \pic(\mathcal A)) \implies \pi_{t-s} \pic(A)^{hG}, \label{notation:Picard sseq} \numberthis \]
and the $E_2$-term can often be identified with continuous cohomology of the $G$-module $\pi_t \pic(A)$.

\item We work at a fixed prime $p$ and height $h$ (mostly one). As such, $p$ and $h$ are often implicit in the notation: e.g. we use the symbol $\mathbb S$ for the both the sphere and its $p$-completion, according to context. To avoid ambiguity, we are explicit in some cases: for example, $KU$ will always mean \emph{integral} K-theory.

\item When indexing spectral sequences, we will always use $s$ for filtration, $t$ for internal degree, and $t-s$ for stem. We abbreviate ``homotopy fixed points spectral sequence'' to HFPSS. We write ``Picard spectral sequence'' for the descent spectral sequence \cref{notation:Picard sseq}.

\item We fix once and for all a regular cardinal $\kappa$ such that $(i)$ $\mathbb S_{K(1)} \in \mathcal Sp_{K(1)}$ is $\kappa$-compact; $(ii)$ $|\G_1| < \kappa$. The Brauer space $\Brr(\mathcal Sp_{K(1)})$ will by definition be the Picard space of $\mathcal Sp_{K(1)} \in \CAlg(\Pr^L_\kappa)$, which is a \emph{small} space since $\Pr^L_\kappa$ is presentable. As noted in \cref{lem:Br = Br'}, for \emph{relative} Brauer classes this is no restriction.
\end{itemize}

\printunsrtglossary[
    title = {\normalsize List of symbols}
]

\section{The descent spectral sequence}
\label{sec:sseq}

In this short section, we record the descent spectral sequence that will be the starting point for our computations. At any characteristic $(p,h)$, this arises as the descent spectral sequence for the sheaf
    \[ \pic(\mathcal E) \in \Sh(\proet \G, \mathcal Sp_{\geq 0}) \]
of \cite[§3.2]{profinitedescent}. The main input from \opcit\ is the following theorem:

\begin{thm}[\cite{profinitedescent}, Theorem~A and Proposition~5.11]
\begin{enumerate}
    \item There is a strongly convergent spectral sequence
        \[ E_2^{s,t} = H^s(\G, \pi_t \pic(\E)) \implies \pi_{t-s} \pic(\mathcal Sp_\K). \label{eqn:picard sseq} \numberthis \]
        
    \item Its $(-1)$-stem converges to $\pi_0 (B\Picc(\E))^{h\G}$, which contains $\Br_h^0$ as a subgroup. 
    
    \item Differentials on the $E_r$ page agree with those in the $\K$-local $\E$-Adams spectral sequence in the region $t \geq r + 1$, and for classes $x \in E_r^{r,r}$ we have
        \[ d_r (x) = d_r^{ASS}(x) + x^2. \]
\end{enumerate}
\end{thm}

In \cite[§4]{profinitedescent}, we used this spectral sequence to recover the computation of $\Pic_1 \coloneqq \Pic(\mathcal Sp_{K(1)})$ (due to \cite{hopkins-mahowald-sadofsky}). In this case, Morava E-theory is the $p$-completed complex K-theory spectrum $KU_p$, acted upon by $\G \cong \Z[p]^\times$ via Adams operations $\psi^a$.

\begin{defn}
We will write $\Brr'(\mathcal Sp_\K \mid \E) \coloneqq (B\Picc(\E))^{h\G}$ for the global sections of the (sheafified) pro\'etale sheaf $\Brr(\mathcal E \mid \E)$, and $\Br'(\mathcal Sp_\K \mid \E) \coloneqq \pi_0 \Brr'(\mathcal Sp_\K \mid \E)$, the group computed in the $(-1)$-stem of \cref{eqn:picard sseq}. Thus $\Br_h^0 \subset \Br'(\mathcal Sp_\K \mid \E)$. At height one, we will show in \cref{sec:azumaya} that these groups agree.
\end{defn}

\subsection{Odd primes}
\begin{figure}
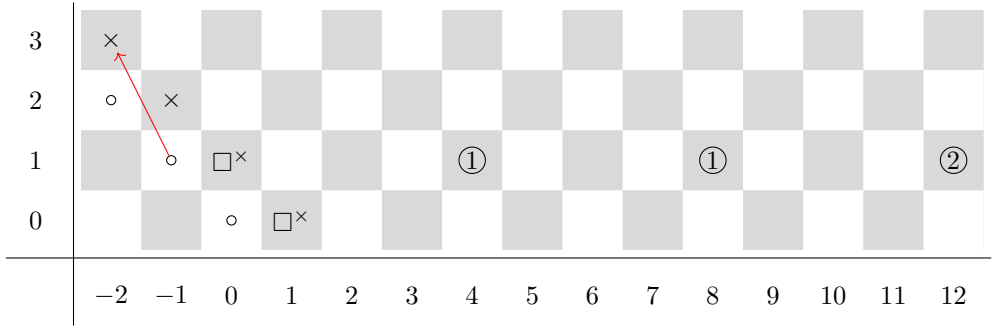

    \centering
    \printpage[grid = chess, xrange = {-2}{12}, yrange={0}{3}, page = 2, name = PicardSSOdd]
    \caption[Picard spectral sequence, odd primes]{The height one Picard spectral sequence for odd primes (implicitly at $p = 3$). Classes are labelled as follows: $\circ = \Z/2$, $\square^\times = \Z[p]^\times$, $\times = \mu_{p-1}$, and circles denote $p$-power torsion (labelled by the torsion degree). Since $\Pic_1 \cong \Pic_1^{\mathrm{alg}} \cong \Z[p]^\times \times \Z/2$, no differentials can hit the $(-1)$-stem. Differentials with source in stem $t-s \leq -2$ have been omitted.}
    \label{fig:Picard spectral sequence p odd}
\end{figure}

We first consider the case $p > 2$. The starting page of the Picard spectral sequence is recorded below:

\begin{lem}[\cite{profinitedescent}, Lemma~4.15]
    At odd primes, the starting page of the descent spectral sequence is given by
    \begin{align}
        E_2^{s,t} = H^s(\Z[p]^\times, \pi_t \pic(KU_p)) = \left\{ \begin{array}{ll}
             \Z/2 & t = 0 \text{ and } s \geq 0 \\
             \Z[p]^\times & t = 1 \text{ and } s = 0,1 \\
             \mu_{p-1} & t = 1 \text{ and } s \geq 2 \\
             \Z/p^{\nu_p(t) + 1} & t = 2(p-1)t'+1 \neq 1 \text{ and } s = 1
        \end{array} \right.
    \end{align} 
    This is displayed in \cref{fig:Picard spectral sequence p odd}. In particular, the spectral sequence collapses for degree reasons at the $E_3$-page.
\end{lem}

\begin{prop} \label{lem:Picard sseq p odd}
At odd primes, $\Br_1^0$ is isomorphic to a subgroup of $\mu_{p-1}$.
\end{prop}

\begin{proof}
The only possible differentials are $d_2$-differentials on classes in the $(-1)$-stem; note that there are no differentials \emph{into} the $(-1)$-stem, since every $E_2$-class in the $0$-stem is a permanent cycle. The generator in $E_2^{1,0}$ supports a $d_2$, since this is the case for the class in $E_2^{1,0}$ of the descent spectral sequence for the $C_2$-action on $KU$ \cite[Prop.~7.15]{gl} (this is displayed in \cref{fig:Brauer spectral sequence for KO}), and the span of Galois extensions
\[\begin{tikzcd}[ampersand replacement=\&] \label{diag:span of Galois extensions}
	{\bm 1_{\K}} \& {KO_p} \& KO \\
	{KU_p} \& {KU_p} \& KU
	\arrow[from=1-2, to=2-2]
	\arrow[from=1-1, to=2-1]
	\arrow[from=1-1, to=1-2]
	\arrow[from=2-1, to=2-2]
	\arrow[from=1-3, to=1-2]
	\arrow[from=2-3, to=2-2]
	\arrow[from=1-3, to=2-3]
\end{tikzcd}\]
allows us to transport this differential (see also \cref{fig:Picard spectral sequence for KO_p to KU_p p odd}). Note that the induced span on $E_2$-pages is
\[\begin{tikzcd}[ampersand replacement=\&]
	{\Z/2} \& {\Z/2} \& {\Z/2} \\
	{\mu_{p-1}} \& {\mu_{p-1}} \& {\Z/2}
	\arrow[from=1-1, to=2-1, red, "d_2"]
	\arrow[from=1-2, to=2-2, red, "d_2"]
	\arrow["{=}"', from=1-3, to=2-3, red, "d_2"]
	\arrow[hook', from=2-3, to=2-2]
	\arrow[no head, from=2-1, to=2-2, equals]
	\arrow[no head, from=1-1, to=1-2, equals]
	\arrow[no head, from=1-3, to=1-2, equals]
\end{tikzcd}\]
in bidegrees $(s,t) = (1,0)$ and $(3,1)$ respectively. Thus
    \[ \Br'(\mathcal Sp_\K \mid \E) \cong \mu_{p-1}. \qedhere \]
\end{proof}

In \cref{sec:cyclic} we will show that this bound is achieved using the cyclic algebra construction of \cite{baker-richter-szymik}; abstractly, this also follows from the detection result of \cref{sec:generators}.

\subsection{The case \texorpdfstring{$p=2$}{p=2}}
We now proceed with the computation of the $(-1)$-stem for the even prime.

\begin{figure}
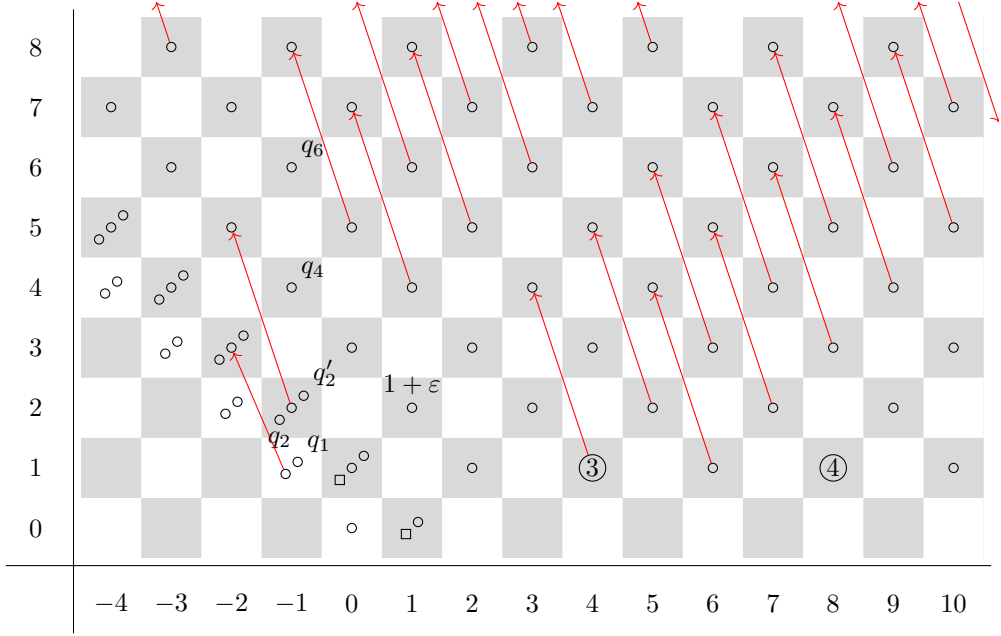

    \centering
    
    \printpage[x range = {-4}{10}, y range = {0}{8}, grid = chess, page = 2--5, name = Brauersseven, classes = {show name = above}]
    
    \caption[Picard spectral sequence, $p=2$]{The Picard spectral sequence for the Galois extension $\bm 1_\K \to \E = KU_2$ at $p=2$. We know that all remaining classes in the $0$-stem survive, by comparing to the algebraic Picard group. Thus the only differentials that remain to compute are those \emph{out} of the $(-1)$-stem; those displayed can be transported from the descent spectral sequence for $\Picc(KO)^{hC_2}$---see \cref{fig:Brauer spectral sequence for KO,fig:Brauer spectral sequence for KO2}. We have not displayed possible differentials out of stem $\leq -2$.}
    
    \label{fig:Brauer spectral sequence p=2}
\end{figure}

\begin{figure}
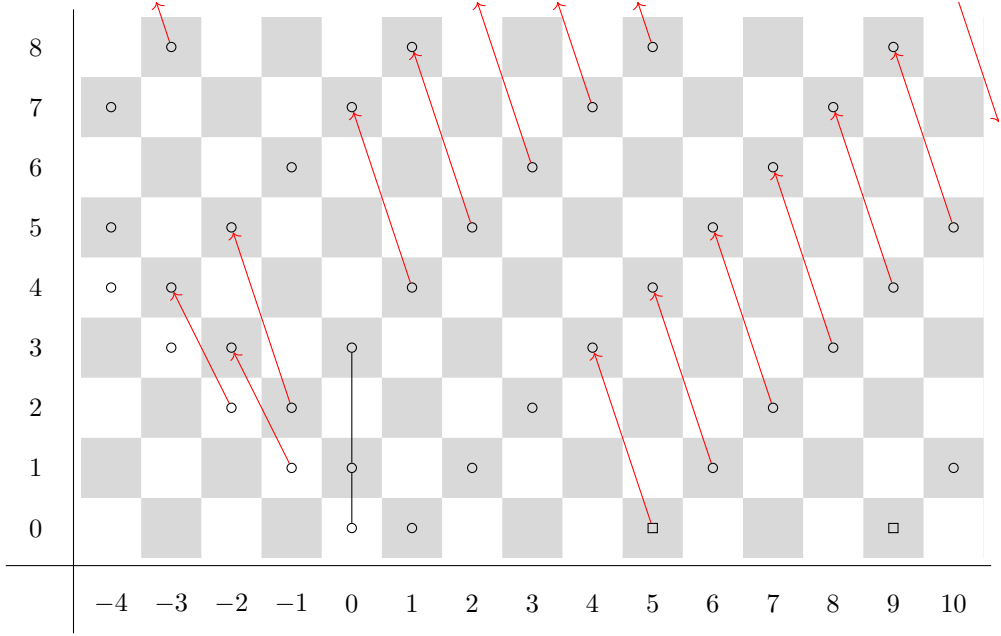

    \centering

    \printpage[x range = {-4}{10}, y range = {0}{8}, grid = chess, page = 2--3, name = BrauerKO]
    
    \caption[HFPSS for $\pic(KO)$]{The $E_3$-page of the Picard spectral sequence for $KO$, as in \cite[Figure~7.2]{gl}.}
    \label{fig:Brauer spectral sequence for KO}
\end{figure}

\begin{lem}[\cite{profinitedescent}, Lemma~4.17]
We have
\begin{align*}
        H^s(\Z[2]^\times, \Pic(KU_2)) &= \left\{ \begin{array}{ll}
          \Z/2 & s = 0 \\
          (\Z/2)^2 & s \geq 1
    \end{array} \right. \\
    H^s(\Z[2]^\times, (\pi_0 KU_2)^\times) &= \left\{ \begin{array}{ll}
          \Z[2] \oplus \Z/2 & s = 0\\
          \Z[2] \oplus (\Z/2)^2 & s = 1 \\
          (\Z/2)^3 & s \geq 2
    \end{array} \right.
\end{align*}
The resulting spectral sequence is displayed in Figure \ref{fig:Brauer spectral sequence p=2}.
\end{lem}

\begin{prop}
At the prime two,
    \[ | \Br_1^0| \leq 32. \]
\end{prop}

\begin{proof}
In \cite[§4]{profinitedescent}, we determined the following differentials:
\begin{itemize}
    \item in degrees $t \geq 3$, differentials agree with the well-known pattern of Adams differentials (e.g. \cite[Figure~3]{beaudry-goerss-henn_splitting}).
    \item the class in bidegree $(s,t) = (3,3)$, which supports a $d_3$ in the Adams spectral sequence, is a permanent cycle.
\end{itemize}
By comparing with the Adams spectral sequence, any classes in the $(-1)$-stem that survive to $E_\infty$ are in filtration at most six; on the $E_2$-page, there are seven such generators. By comparing to the HFPSS for $\Brr'(KO_2 \mid KU_2) = (B\Picc(KU_2))^{hC_2}$ as in \cref{diag:span of Galois extensions}, we obtain the following differentials:
\begin{itemize}
    \item a $d_2$ on the class in $H^1(C_2, \Pic(KU_2)) \subset H^1(\Z[2]^\times, \Pic(KU_2))$,
    \item a $d_3$ on the class in $H^2(C_2, (\pi_0 KU_2)^{\times}) \subset H^2(\Z[2]^\times, (\pi_0 KU_2)^{\times})$.
\end{itemize}
This gives the claimed upper bound.
\end{proof}

In \cref{sec:azumaya} we will show that this bound is \emph{also} achieved.

\begin{remark} \label{rem:E2 generators}
For later reference, we name the following generators:

\begin{enumerate}
    \item $q_1 \in E_2^{1,0}$ is the generator of $H^1(1 + 4\Z[2], \Z/2) \subset H^1(\Z[2]^\times, \Z/2)$,

    \item $q_2 \in E_2^{2,1}$ is the generator of $H^2(C_2, 1 + 4\Z[2]) \subset H^2(\Z[2]^\times, \Z[2]^\times)$,

    \item $q_2' \in E_2^{2,1}$ is the generator of $H^1(C_2, \Z[2]^\times) \otimes H^1(\Z[2], \Z[2]^\times) \subset H^2(\Z[2]^\times, \Z[2]^\times)$,

    \item $q_4$ is the unique class in $E_2^{4,3}$,

    \item $q_6$ is the unique class in $E_2^{6,5}$.
\end{enumerate}
\end{remark}

While $q_6$ survives to $E_\infty$ by sparsity in \cref{fig:Brauer spectral sequence p=2}, the other classes are sources of possible differentials. We will show that in fact all are permanent cycles.

\begin{figure}
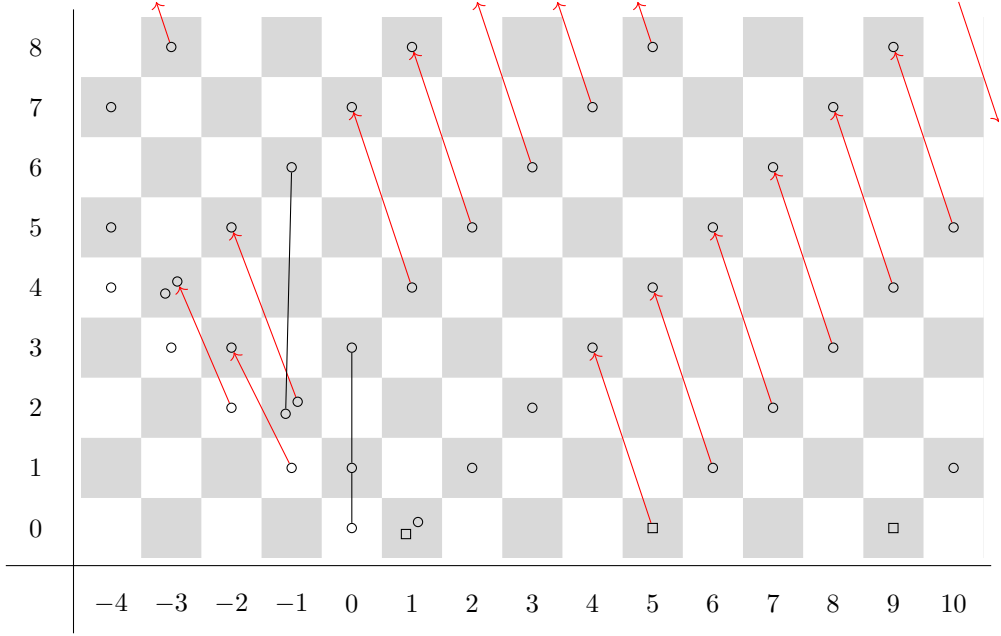

    \centering

    \printpage[x range = {-4}{10}, y range = {0}{8}, grid = chess, page = 2--3, name = BrauerKO2]
    
    \caption[HFPSS for $\pic(KO_2)$]{The Picard spectral sequence for $KO_2$. Differentials come from the comparison with $KO$, and the resulting map of spectral sequences. See \cref{sec:azumaya} for the extension in the $(-1)$-stem.}
    \label{fig:Brauer spectral sequence for KO2}
\end{figure}

\section{Descent for compact generators}
\label{sec:generators}

Given a Galois extension $\bm 1 \to A$ in a stable homotopy theory $\mathcal C$, we showed in \cite[§5]{profinitedescent} that the Picard spectral sequence computes the subgroup
    \[ \Br'(\bm 1 \mid A) \subset \Br'(A) = \Pic(\Cat_A) \]
of the cohomological Brauer group. To relate this to the Brauer-Azumaya group classifying Azumaya algebras in $\mathcal C$, we prove a descent result for compact generators valid in the $\K$-local setting. This is entirely analogous to the theory of \cite[§6.3]{antieau-gepner} and \cite[§6.4]{gl}.

\begin{defn}
Let $\mathcal C$ be a presentably symmetric monoidal \inftycat{}, and $\mathcal D \in \Mod_{\mathcal C}(\Pr^L)$. An object $D \in \mathcal D$ is \emph{$\mathcal C$-compact} if the functor
    \[ \eMap_{\mathcal D}(D, - ) : \mathcal D \to \mathcal C \numberthis \label{eqn:enriched mapping functor} \]
preserves filtered colimits. $D$ is an \emph{$\mathcal C$-compact generator} if it is $\mathcal C$-compact. We say that $D$ is a \emph{$\mathcal C$-generator} if the functor \cref{eqn:enriched mapping functor} is conservative; when $\mathcal C$ is stable, it is equivalent that \cref{eqn:enriched mapping functor} detects zero objects. A \emph{$\mathcal C$-compact generator} of $\mathcal D$ is an object $D \in \mathcal D$ that is both $\mathcal C$-compact and a $\mathcal C$-generator, and we shall write $\mathcal D\cg \subset \mathcal D$ for the full subcategory of such $D$ \footnote{These might also be called \emph{e}nriched \emph{c}ompact \emph{g}enerators.}.
\end{defn}

\begin{ex}
The $\K$-local sphere is an $\mathcal Sp_{\K}$-compact generator of $\mathcal Sp_{\K}$. More generally, we always have $\bm 1 \in \mathcal C\cg$ since in this case \cref{eqn:enriched mapping functor} is the identity functor.
\end{ex}

Our first objective is to show that Schwede-Shipley theory goes through in the presence of a $\mathcal C$-compact generator.

\begin{defn}
Let $\mathcal C$ be a stable homotopy theory. We say $\mathcal C$ is \emph{rigidly generated} if it is generated under colimits by dualisable objects. That is, the localising category generated by $\mathcal C^{\mathrm{dbl}}$ is $\mathcal C$ itself.
\end{defn}

\begin{ex}
\begin{enumerate}
    \item $\mathcal Sp$ is generated under colimits by shifts of $\bm 1$, and so rigidly generated.
    
    \item If $\mathcal C$ is a rigidly generated stable homotopy theory and $L \colon \mathcal C \to \mathcal C'$ a monoidal localisation, then $\mathcal C'$ is rigidly generated. Thus $\mathcal Sp_{\K}$ is rigidly generated.
    
    \item For a compact lie group $G$, the \inftycat{} $\mathcal S^G$ of $G$-spaces is generated under colimits by orbits $G/H$ (e.g. \cite[Theorem~1.8]{mandell-may}). Its stabilisation $\mathcal Sp_{\mathscr U}^G$ at any $G$-universe $\mathscr U$ (as defined in \cite[Corollary~C.7]{gepner-meier}) is generated under colimits by shifts $\Sigma^{-V} \Sigma^\infty_{\mathscr U} G/H_+$ as $V$ ranges over representations in $\mathscr U$, and if $\mathscr U$ is complete then these are invertible by virtue of the Wirthm\"uller isomorphism \cite[(4.16)]{greenlees-may_equivariant}. Thus the \inftycat{} $\mathcal Sp^G$ of genuine $G$-spectra is rigidly generated.
\end{enumerate}
\end{ex}

\begin{prop}[Enriched Schwede-Shipley] \label{lem:enriched barr-beck}
Let $\mathcal C$ be a rigidly generated stable homotopy theory and $\mathcal D \in \Cat_{\mathcal C}$. Suppose that $D \in \mathcal D\cg$, and write $A \coloneqq \underline{\End}_{\mathcal D}(D) \in \Alg(\mathcal C)$. Then there is an $\mathcal C$-linear equivalence
    \[ \mathcal D \simeq \LMod_{A}(\mathcal C). \]
\end{prop}

\begin{proof}
The object $D \in \mathcal D$ determines canonically a $\mathcal C$-linear left adjoint $F \colon \mathcal C \to \mathcal D$, with right adjoint $G \coloneqq \eMap_{\mathcal D}(D, -)$. According to \cite[Proposition~4.8.5.8]{ha}, it is enough to check the following:
\begin{enumerate}
    \item $G$ preserves colimits of simplicial objects: in fact $G$ preserves all colimits. Indeed, $G$ preserves filtered colimits since $D$ is $\mathcal C$-compact, and finite colimits as it is a right adjoint. 
    
    \item $G$ is conservative: this is by definition of $\mathcal C$-compact generators.
    
    \item for every $D' \in \mathcal D$ and $C \in \mathcal C$, the map
        \[ C \otimes FG(D') = C \otimes G(D') \otimes D \to C \otimes D' \]
    is adjoint to an equivalence
        \[ C \otimes G(D') \xrightarrow{\sim} G(C \otimes D'). \numberthis \label{eqn:projection formula for enriched Schwede-Shipley} \]
    But by $(i)$, the functor $G$ preserves all colimits, so by rigid generation we reduce to $C$ dualisable. In this case, \cref{eqn:projection formula for enriched Schwede-Shipley} is the composite equivalence
    \begin{align*}
        C \otimes \eMap_{\mathcal D}(D,D') &\simeq \eMap_{\mathcal C}(C^\vee, \eMap_{\mathcal D}(D,D')) \\
        &\simeq \eMap_{\mathcal D}(C^\vee \otimes D,D') \\
        &\simeq \eMap_{\mathcal D}(D, C \otimes D'). \qedhere
    \end{align*}
\end{enumerate}

\end{proof}

To use this to produce Azumaya algebras, we first need to be able to produce $\mathcal C$-compact generators. For this we will prove the following result:

\begin{prop}[Descent for $\mathcal C$-compact generators] \label{lem:descent for compact generators}
Let $\bm 1 \to A$ be a descendable extension in a rigidly generated stable homotopy theory $\mathcal C$, and suppose that $A$ is dualisable. Then $\mathcal D \in \Cat_{\mathcal C}$ admits a $\mathcal C$-compact generator if and only if $\Mod_A(\mathcal D) \coloneqq \mathcal D \otimes_{\mathcal C} \Mod_A(\mathcal C)$ admits one.
\end{prop}

The proof relies on the following basic lemma:

\begin{lem}
Suppose $\mathcal C \in \CAlg(\Pr^L)$ and $A \in \mathcal C$ is dualisable, then $A$ is faithful if and only if $A^\vee$ is.
\end{lem}

\begin{proof}
Assume that $A$ is faithful, and that $A^\vee \otimes X = 0$; the converse is given by taking duals. Then the identity on $A \otimes X$ factors as
    \[ A \otimes X \to A^\vee \otimes A \otimes A^\vee \otimes X \to A \otimes X, \]
and in particular $A \otimes X = 0$. By faithfulness of $A$, this implies $X = 0$.
\end{proof}

\begin{proof}[Proof (\cref{lem:descent for compact generators}).]
As in \cite{gl}, we will make use of the functors
\[\begin{tikzcd}[ampersand replacement=\&]
	{\mathcal C} \& {\Mod_A(\mathcal C)}
	\arrow["{i^*}"{description}, from=1-2, to=1-1]
	\arrow["{i_!}", shift left=3, from=1-1, to=1-2]
	\arrow["{i_*}"', shift right=3, from=1-1, to=1-2]
\end{tikzcd}\]
and the adjunctions (denoted by the same symbols) between $\mathcal D$ and $\Mod_A(\mathcal D)$. In fact, we claim these adjunctions are $\mathcal C$-linear.

Assuming this, we give the proof of the proposition. If $\mathcal D$ admits a $\mathcal C$-compact generator $D$, it is straightforward to check that $i_! D \in \Mod_A(\mathcal D)\cg$: indeed, $i_!D$ is $\mathcal C$-compact because $D$ is so and $i^*$ preserves colimits, while $i_!D$ generates because $D$ does and $i^*$ is conservative. Conversely, suppose we have $D \in \Mod_A(\mathcal D)\cg$, and consider $i^* D \in \mathcal D$. By dualisability of $A$, the right adjoint
    \[ i_* = \eMap_{\mathcal C}(A, -) \simeq A^\vee \otimes - : \mathcal C \to \Mod_A(\mathcal C) \]
preserves colimits, and hence the right adjoint $i_* \colon \Mod_A(\mathcal D) \to \mathcal D$ does too. As a result, $i^* D$ is $\mathcal C$-compact. On the other hand, if $X \in \mathcal D$ and $\eMap_{\mathcal D}(i^* D, X) = 0$, then $\eMap_{\Mod_A(\mathcal D)}(D, i_* X) = 0$ and so
    \[ i_* X = A^\vee \otimes X = 0. \]
Now faithfulness of $A^\vee$ implies that $X = 0$.

It remains to prove that the adjunctions $i_! \dashv i^* \dashv i_*$ are $\mathcal C$-linear.  For $i_!$ we observed $\mathcal C$-linearity in the proof of \cref{lem:enriched barr-beck}. To see $\mathcal C$-linearity for $i^*$ it is enough to prove that the canonical map
    \[ \theta : C \otimes i^* M \to i^*(C \otimes M) \label{eqn:canonical map for show p^* is C-linear} \]
is an equivalence for every $C$ and $M$, and by rigid generation we reduce to $C$ dualisable. As in \cite[Remark~D.7.4.4]{sag} one checks that $\Map_{\mathcal C}(C', \theta)$ is the composite equivalence
\begin{align*}
    \Map(C', C \otimes i^* M) &\simeq \Map(C' \otimes C^\vee, i^* M)  \\
    &\simeq \Map(i_!(C' \otimes C^\vee), M) \\
    &\simeq \Map(i_!C' \otimes C^\vee, M) \\
    &\simeq \Map(i_!(C'), C \otimes M) \\
    &\simeq \Map(C', i^*(C \otimes M))
\end{align*}
for any $C' \in \mathcal C$, which gives the claim.
\end{proof}

\begin{ex}
If $A \to B$ is an $E$-local Galois extension of ring spectra with stably dualisable Galois group $G$, then Rognes \cite[Proposition~6.2.1]{rognes} shows that $B$ is dualisable over A. For example, this covers the following cases:
\begin{enumerate}
    \item $E = \mathbb S$ and $G$ is finite or compact Lie.
    \item $E = \F[p]$ and $G$ is $p$-compact.
    \item $E = \K$ and $G = K(\pi, m)$ for $\pi$ a finite $p$-group and $m \leq h$.
\end{enumerate}
\end{ex}

\begin{cor}[$\Br = \Br'$] \label{lem:Br = Br'}
Let $\bm 1 \to A$ be a faithful dualisable Galois extension in a rigidly generated stable homotopy theory $\mathcal C$, and $\mathcal Q \in \pi_0 B\Picc(A)^{hG}$ a relative Brauer-Grothendieck class. Then $\mathcal Q$ is represented by some Azumaya algebra $Q$ whose basechange to $A$ is (Morita) trivial: that is,
    \[ \mathcal Q \simeq \Mod_Q(\mathcal C) \in \Mod_{\mathcal C}({\Pr}^L), \]
and $\Mod_{A \otimes Q}(\mathcal C) \simeq \Mod_A(\mathcal C)$. Thus the map
    \[ \Br(\bm 1 \mid A) \to \Br'(\bm 1 \mid A) \]
is an isomorphism.
\end{cor}

\begin{proof}
We claim that $\mathcal Q$ is $\kappa$-compactly generated, so that $\mathcal Q \in \Pic(\Cat_{\mathcal C})$. Given this, the result follows from \cref{lem:descent for compact generators}: by assumption, $\Mod_A(\mathcal Q) \simeq \Mod_A(\mathcal C)$, and so
    \[ A \in \Mod_A(\mathcal C)\cg \simeq \Mod_A(\mathcal Q)\cg. \]
By descent for compact generators we obtain $D \in \mathcal Q\cg$, and so Schwede-Shipley theory yields a $\mathcal C$-linear equivalence $\mathcal Q \simeq \Mod_Q(\mathcal C)$, where $Q = \underline{\End}_{\mathcal Q}(D)$. 

For the claim, note that as in \cite[Corollary~3.42]{mathew_galois}, $\mathcal Q$ is the limit of the cosimplicial diagam
    \[ \cosimp[\mathcal Q \otimes_{\mathcal C} \Mod_A(\mathcal C)]{\mathcal Q \otimes_{\mathcal C} \Mod_{A^{\otimes 2}}(\mathcal C)}, \]
and for $j \geq 1$ we have $\mathcal Q \otimes \Mod_{A^{\otimes j}}(\mathcal C) \simeq \Mod_{A^{\otimes j}} (\mathcal C) \in \Pr^L_{\kappa}$. Thus $\mathcal Q \in \Pr^L_{\kappa}$.
\end{proof}

\begin{remark}
In a previous version, we claimed that Morava E-theory is dualisable in $\mathcal Sp_\K$. As pointed out to us by Maxime Ramzi, while $\E$ is Spanier-Whitehead self-dual \cite{strickland_gross-hopkins, beaudry-goerss-hopkins-stojanoska} and hence \emph{reflexive}, it is not dualisable: for example, $\K_* \E$ would otherwise be finite by \cite[Theorem~8.6]{hovey-strickland}. We will bypass this issue at height one by showing that all generators of $\Br'(\mathcal Sp_\K \mid \E)$ are in fact trivialised in a \emph{finite} Galois extension of the sphere, and hence lift to Azumaya algebras by \cref{lem:Br = Br'}; we do not know if $\Br_h^0 \cong \Br'(\mathcal Sp_\K \mid \E)$ at arbitrary height.
\end{remark}

\section{Explicit generators}
\label{sec:cyclic}

In this section, we give some explicit constructions of Azumaya algebras from Galois extensions. Most of this section works in an arbitrary stable homotopy theory $\mathcal C$. We will use these constructions in \cref{sec:azumaya} to describe generators of the group $\Br_1^0$, and hence to solve extension problems.

\subsection{\texorpdfstring{$\Z[p]$}{Zp}-extensions}
We will begin with a straightforward construction for extensions with Galois group $\Z[p]$, using the fact that $\cd(\Z[p]) = 1$.

\begin{remark}
Let $G$ be a profinite group, and $\mathcal B \in \Sh(\proet G, \mathcal S_*)$. Under suitable assumptions on $\mathcal B$, d\'ecalage gives an isomorphism between the descent spectral sequence and the spectral sequence for the \v Cech nerve of $G \to *$ \cite[Appendix~A]{profinitedescent}. Moreover, the homotopy groups of the latter can often be identified with the complex of continuous cochains with coefficients in $\pi_t B$, yielding an isomorphism on the $E_2$-pages
    \[ H^*(G, \pi_t B) \to H^*(\proet G, \pi_t \mathcal B). \label{eqn:group cohomology to proetale cohomology} \numberthis \]
For example, this is the case for the sheaf $\pic(\mathcal E)$ for any Morava $E$-theory $E(k, \Gamma)$.
\end{remark}

\begin{lem}
Let $G = \Z[p]$ or $\widehat{\Z}$, and write $\zeta \in G$ for a topological generator. Suppose $\mathcal C$ is a stable homotopy theory and $\mathcal B \in \Sh(\proet{G}, \mathcal S_*)$, with $B \coloneqq \mathcal B(G/*)$ and $B^{hG} \coloneqq \Gamma \mathcal B$. If the canonical map \cref{eqn:group cohomology to proetale cohomology} is an isomorphism, then 
        \[ B^{hG} \simeq \operatorname{Eq} \left(id, \zeta: B \rightrightarrows B \right). \label{eqn:Zp fixed points is equaliser} \numberthis \]
\end{lem}

\begin{proof}
Write $B'$ for the equaliser in \cref{eqn:Zp fixed points is equaliser}. The $G$-map $(id, \zeta) \colon G \to G \times G$ gives rise to maps
    \[ B \rightrightarrows \mathcal B(G \times G) \to B \label{eqn:Zp fixed points to equaliser} \numberthis \]
factoring $id, \zeta \colon B \rightrightarrows B$, and the identification $B^{hG} \simeq \lim \mathcal B(G^{\bullet + 1})$ gives a distinguished nullhomotopy in \cref{eqn:Zp fixed points to equaliser} after precomposing with the coaugmentation $\eta \colon B^{hG} \to B$. Thus $\eta$ factors through $\theta \colon B^{hG} \to B'$. Taking fibres, the descent spectral sequence for $\mathcal B$ implies that $\pi_* \theta$ fits in a commutative diagram
\[\begin{tikzcd}[column sep=small]
	0 & {(\pi_{t+1}B)_G} & {\pi_tB^{hG}} & {(\pi_tB)^G} & 0 \\
	0 & {(\pi_{t+1}B)_G} & {\pi_t B'} & {(\pi_tB)^G} & 0
	\arrow[from=1-1, to=1-2]
	\arrow[from=1-2, to=1-3]
	\arrow[from=1-3, to=1-4]
	\arrow[from=1-4, to=1-5]
	\arrow[from=2-1, to=2-2]
	\arrow[from=2-2, to=2-3]
	\arrow[from=2-3, to=2-4]
	\arrow[from=2-4, to=2-5]
	\arrow[no head, equals, from=1-2, to=2-2]
	\arrow[no head, equals, from=1-4, to=2-4]
	\arrow["\pi_t \theta", from=1-3, to=2-3]
\end{tikzcd} \]
and is therefore an equivalence.
\end{proof}

\begin{constr} \label{lem:Azumaya algebras from H1 classes}
Suppose $\bm 1 \to A$ is a descendable Galois extension in a stable homotopy theory $\mathcal C$, with group $G = \Z[p]$ or $\widehat{\Z}$. Then $\mathcal B \coloneqq B\Picc(A) \in \Sh(\proet G, \mathcal S_*)$, and in good cases $\mathcal B$ satisfies the assumption that \cref{eqn:group cohomology to proetale cohomology} be an isomorphism: for example, this is the case whenever each $\pi_t \pic(A)$ is the limit of a tower of finite sets by \cite[Lemma~4.3.9]{bs}. Thus
    \[ \Brr'(\bm 1 \mid A) \simeq \mathrm{Eq} \left(id, \zeta^* :B\Picc(A) \rightrightarrows B\Picc(A) \right), \]
and so a relative cohomological Brauer class is given by the data of an $A$-linear equivalence
    \[ \xi: \Mod_A(\mathcal C) \xrightarrow{\sim} \zeta^* \Mod_A(\mathcal C). \]
In fact, if $X \in \Pic(A)$, then we may form the $A$-linear composite
    \[ \zeta_!^X : \Mod_A \xrightarrow{X \otimes_A -} \Mod_A \xrightarrow{\zeta_!} \zeta^* \Mod_A. \]
This gives an isomorphism
    \[ (A,-) \coloneqq \zeta_!^{(-)} : \Pic(A)_G \cong H^1(G, \Pic(A)) \to \Br'(\bm 1 \mid A). \]
\end{constr}

\begin{ex} \label{lem:KO2-trivial algebra from H1 cocycle}
At the prime $2$, the extension $\bm 1_\K \to KO_2$ is a descendable $\Z[2]$-Galois extension. As a result, \cref{lem:Azumaya algebras from H1 classes} applies, and we can form the Brauer class $(KO_2, X)$ associated to any $X \in \Pic(KO_2)$ as above. For example, since $KO_2$ is $8$-periodic one can form the cohomological Brauer class
    \[ (KO_2, \Sigma KO_2) 
    \in \Br'(\bm 1_\K \mid KO_2). \]
\end{ex}

\begin{ex} \label{lem:KO2 algebra from H1 cocycle}
Let $KO_2\un \coloneqq \varinjlim_n (KO_2)_{W(\F[2^n])}$ be the ind-\'etale $KO_2$-algebra given by the maximal unramified extension of $\pi_0 KO_2 = \Z[2]$; since \'etale extensions are uniquely determined by their $\pi_0$, one can also describe this as
    \[ KO_2\un = KO_2 \otimes_{\mathbb S} \mathbb{SW}, \]
where $\mathbb{SW} = W^+(\overline{\F}_2)$ denotes the spherical Witt vectors \cite[§5.2]{ellipticii}. The extension $KO_2 \to KO_2\un$ is a descendable Galois extension: indeed, $KO_2 \otimes_{\mathbb S} (-) $ preserves finite limits, and $\mathbb S \to \mathbb{SW}$ is descendably $\widehat{\Z}$-Galois. As a result, \cref{lem:Azumaya algebras from H1 classes} applies, and we can form the Brauer class $(KO_2 \un, X)$ associated to any $X \in \Pic(KO_2\un)$ as above. For example, since $KO_2\un$ is $8$-periodic one can form the cohomological Brauer class
    \[ (KO_2\un, \Sigma KO_2\un) 
    \in \Br'(KO_2 \mid KO_2\un). \]
In fact, $(KO_2\un, \Sigma KO_2\un)$ is an element of the \emph{\'etale locally trivial} Brauer group $\mathrm{LBr}(KO_2)$ of \cite{antieau-meier-stojanoska}; we discuss this in \cref{sec:azumaya}.
\end{ex}

\subsection{Cyclic algebras}
Suppose that $\mathcal C$ is a stable homotopy theory and $\bm 1 \to A$ a finite Galois extension in $\mathcal C$ with group $G$. Suppose also given the following data:
\begin{enumerate}
    \item an isomorphism $\chi: G \cong \Z/k$.
    \item a \emph{strict} unit $u \in \pi_0 \mathbb G_m(\bm 1)$.
\end{enumerate}
In this section, we will use this to define a relative Azumaya algebra $(A, \chi, u) \in \Br(\mathcal C \mid A)$.

Let us first recall the construction when $\mathcal C$ is the category of modules over a classical ring. Then we begin with a $G$-Galois extension $R \to A$ of rings, and define a $G$-action on the matrix algebra $M_k(A)$ as follows: we let $\sigma \coloneqq \chi^{-1}(1)$ act as conjugation by the matrix
\begin{align} \label{eqn:widetilde u in GLk R}
    \widetilde{u} \coloneqq \left[ \begin{array}{cccc}
         0 & & & u \\
         1 & 0 & & \\
         & \ddots & \ddots & \\
         && 1 & 0
    \end{array} \right]
\end{align}

Since $\widetilde{u}^k = u I_k \in \GL_k A$ is central, this gives a well-defined action on $M_k A$. We can use this to twist the usual action of $G$ on the matrix algebra; passing to fixed points, we obtain the cyclic algebra
\begin{align} \label{eqn:defn of cyclic algebras}
    (A, \chi, u) \in \Br(R \mid A).
\end{align}

The construction for general $\mathcal C$ will be directly analogous. Moreover, when $R$ is a field it is well-known (see for example \cite{cs_brauer}) that under the isomorphism
    \[ \Br(R \mid A) \cong H^2(G,A^\times), \]
the cyclic algebra $(A, \chi, u)$ maps to the cup-product $\beta(\chi) \cup u$, where $\beta$ denotes the Bockstein homomorphism
    \[ \chi \in \Hom(G, \Z/k) = H^1(G, \Z/k) \xrightarrow{\beta} H^2(G, \Z), \]
and we use the $\Z$-module structure on $A^\times$. As a result, we obtain an isomorphism
    \[ \widehat{H}^0(G, A^\times) = A^\times/N_e^G A^\times \to \Br(R \mid A) \]
sending $u \mapsto (A, \chi, u)$. We will prove an analogous result for cyclic algebras in arbitrary stable homotopy theories, which will allow us to detect permanent cycles in the descent spectral sequence; conversely, this will allow us to assert that the cyclic algebras we construct are nontrivial.

We now begin the construction of $(A, \chi, u)$ for arbitrary $\mathcal C$.

\begin{defn}
Let $\mathcal C$ be a stable homotopy theory. Given $A \in \CAlg(\mathcal C)$ and $k \geq 1$, we will write
    \begin{itemize}
        \item $M_k(A) \coloneqq \End_{\Mod_{A}(\mathcal C)}(A^{\oplus k})$,
        \item $\GL_k(A) \coloneqq \Aut_{\Mod_A(\mathcal C)}(A^{\oplus k})$,
        \item $\PGL_k(A) \coloneqq \Aut_{\Alg_{A}(\mathcal C)}(M_k(A))$.
    \end{itemize}
These are all $\mathbb E_1$-monoids under composition, and we have $\mathbb E_1$ maps
    \[ \GL_k(A) \to M_k(A) \qquad \text{and} \qquad \GL_k(A) \to \PGL_k(A), \]
the first of which is by definition. The second will be defined immediately below.
\end{defn}

\begin{remark}
If $R$ is an \Einfty-ring, then Gepner and Lawson \cite[Corollary~5.19]{gl} show there is a fibre sequence
\begin{align} \label{eqn:gl fibre sequence for BAut}
    \coprod_{\pi_0 \Mod_R^{\mathrm{cg}}} B\Aut_{A}(M) \to \coprod_{\pi_0 \Azz_R^{\mathrm{triv}}} B\Aut_{\Alg_R}(A) \to B \Picc(R).
\end{align}
Here $\Mod_R^{\mathrm{cg}} \subset \iota \Mod_R$ denotes the classifying space of \emph{compact generators}, and $\Azz_R^{\mathrm{triv}} \subset \iota \Alg_R$ the space of \emph{Morita trivial} Azumaya algebras. To see this, we consider the components over $\Mod_R \in \Br(R)$ in the commutative square
\[\begin{tikzcd}[ampersand replacement=\&]
	{_A \BMod_A^{\mathrm{cg}}} \& {\{\LMod_A\}} \\
	{\Azz_R} \& {\Brr_R}
	\arrow[from=2-1, to=2-2]
	\arrow[from=1-2, to=2-2]
	\arrow[from=1-1, to=2-1]
	\arrow[from=1-1, to=1-2]
\end{tikzcd}\]
noting that this diagram is Cartesian \cite[Proposition~5.17]{gl}: this follows by identifying those objects in
    \[ \Fun_{R}(\LMod_A, \LMod_B) \simeq {}_A \BMod_B \]
which correspond to the $R$-linear equivalences. By \cite[Corollary~2.1.3]{hl}, $\mathcal C$-linear equivalences between module categories in an arbitrary $\mathcal C$ correspond to bimodules $M$ that are \emph{full} and \emph{dualisable}, and so in this context there is an analogous pullback diagram
\[\begin{tikzcd}[ampersand replacement=\&]
	{_A \BMod_A^{\mathrm{fd}}(\mathcal C)} \& {\{\LMod_A (\mathcal C) \}} \\
	{\Azz (\mathcal C)} \& {\Brr (\mathcal C)}
	\arrow[from=2-1, to=2-2]
	\arrow[from=1-2, to=2-2]
	\arrow[from=1-1, to=2-1]
	\arrow[from=1-1, to=1-2]
\end{tikzcd}\]
for any Azumaya algebra $A$. Restricting again to the unit component in $\Brr(\mathcal C)$, we obtain a fibre sequence
    \[ \coprod_{M \in \pi_0 \mathcal C^{\mathrm{fd}}} B\Aut_{A}(M) \to \coprod_{A \in \pi_0 \Azz(\mathcal C)^{\mathrm{triv}}} B\Aut_{\Alg(\mathcal C)}(A) \to B \Picc(\mathcal C). \]

In particular, taking $A = M_k(\bm 1)$ yields the fibre sequence
\begin{align} \label{eqn:fibre sequence analogue of BPGL}
    \coprod_{\End(M) = A} B\Aut_{A}(M) \to B\PGL_k(\bm 1) \to B \Picc(\mathcal C),
\end{align}
and restricting to $M = \bm 1^{\oplus k} \in \mathcal C^{\mathrm{fd}}$ gives the desired map $\GL_k(\bm 1) \to \PGL_k(\bm 1)$ after taking loops.
\end{remark}

\begin{remark}
Suppose that $\mathcal C$ is a stable homotopy theory and $u \in \pi_0 \GL_1(\bm 1) = \pi_1 \Pic(\mathcal C)$. We define a map $\hat{u}: \Z/k \to \PGL_k(\bm 1)$ by virtue of the following commutative diagram, whose bottom row is a shift of the fibre sequence \cref{eqn:fibre sequence analogue of BPGL}:
\[\begin{tikzcd}[ampersand replacement=\&, column sep = small]
	{\Z/k} \& {S^1} \& {S^1} \\
	{\PGL_k(\bm 1)} \& {\Picc(\mathcal C)} \& {B\GL_k(\bm 1)}
	\arrow[from=2-1, to=2-2]
	\arrow[from=2-2, to=2-3]
	\arrow["B\widetilde{u}", from=1-3, to=2-3]
	\arrow["Bu", from=1-2, to=2-2]
	\arrow["k", from=1-2, to=1-3]
	\arrow[from=1-1, to=1-2]
	\arrow[dashed, "\widehat{u}"', from=1-1, to=2-1]
\end{tikzcd} \numberthis \label{diag:definition of widehat u as a plain map} \]
The map
    \[ \widetilde{u} \in \pi_0 \Map_{\mathbb E_1}(\Z, \GL_k (\bm 1)) \cong \pi_0 \GL_k(\bm 1) \cong \GL_k([\bm 1, \bm 1]) \]
is given by the matrix in \cref{eqn:widetilde u in GLk R}---commutativity of the right-hand square above is implied by the computation $\widetilde{u}^k = u I_k \in \GL_k([\bm 1, \bm 1])$. In fact, when the chosen unit is \emph{strict}, we will prove at the end of the section that $\widehat u$ deloops:

\begin{prop} \label{lem:widehat u is E1}
If $u \in \pi_0 \G_m(A)$, the map $\widehat{u}$ lifts naturally to an $\mathbb E_1$ map $\widehat u \colon \Z/k \to \PGL_k(A)$. More precisely, there is a natural transformation $\widehat u \colon \G_m \to \Map_{\mathbb E_1}(\Z/k, \PGL_k)$ which lifts
    \[ \widehat u \in \pi_0 \Map(\Z/k, \PGL_k(\bm 1 [u^{\pm 1}])) \cong \pi_0 \Map_{\mathcal P(\CAlg(\mathcal C)\op)}(\G_m, \Map(\Z/k, \PGL_k)). \]
\end{prop}
\end{remark}

\begin{remark}
The choice of lift $\widehat u$ is only determined up to homotopy in $\Map_{\mathcal P(\CAlg(\mathcal C)\op)}(\G_m, \Map_{\mathbb E_1}(\Z/k, \PGL_k))$.
\end{remark}

\begin{constr}
Suppose that $\bm 1 \to A$ is a finite $G$-Galois extension in $\mathcal C$, and $u \in \pi_0 \G_m(\bm 1)$. We obtain a commutative square of spaces with $G$-action
\[\begin{tikzcd} \label{diag:square of spaces for cyclic algebra formula}
	{B \Z/k} & {B^2\Z} \\
	{B\PGL_k(A)} & {B\Picc(A)}
	\arrow["{B\widehat{u}}"', from=1-1, to=2-1]
	\arrow["\delta"', from=2-1, to=2-2]
	\arrow["B^2 u", from=1-2, to=2-2]
	\arrow["\beta", from=1-1, to=1-2]
\end{tikzcd} \numberthis \]
where the action on the top row is trivial. On homotopy fixed points, we obtain the square
\[\begin{tikzcd} \label{diag:square of fixed point spaces for cyclic algebra formula}
	{B \Z/k^{BG}} & {B^2\Z^{BG}} \\
	{(B\PGL_k(A))^{hG}} & {(B\Picc(A))^{hG}}
	\arrow[from=1-1, to=2-1]
	\arrow[from=2-1, to=2-2]
	\arrow[from=1-2, to=2-2]
	\arrow[from=1-1, to=1-2]
\end{tikzcd} \numberthis \]
whose bottom right term is the relative Brauer space $\Brr(\bm 1 \mid A)$.
\end{constr}

\begin{defn}
Given a Galois extension $\bm 1 \to A$ with group $G$ in a stable homotopy theory $\mathcal C$, and given $\chi: G \cong \Z/k$ and $u \in \pi_0 \G_m(\bm 1)$, define the relative Brauer class $(A, \chi, u) \in \Br(\bm 1 \mid A)$ to be the image of $\chi$ under the composite map in \cref{diag:square of fixed point spaces for cyclic algebra formula},
    \[ \chi \in \Hom(G, \Z/k) = \pi_0 B\Z/k^{BG} \to \pi_0 (B\Picc(A))^{hG} \cong \Br(\bm 1 \mid A). \]
The final isomorphism is inverse to the map $\Br(\bm 1 \mid A) \xrightarrow{\sim} \Br'(\bm 1 \mid A)$ considered in \cref{sec:generators}.
\end{defn}

\begin{prop}
The Brauer class $(A, \chi, u)$ agrees with the class of the cyclic algebra $A(A, \chi, u)$ defined in \cite[§4]{baker-richter-szymik}.
\end{prop}

\begin{proof}
Going down and left in \cref{diag:square of fixed point spaces for cyclic algebra formula}, we see that $(A,\chi,u)$ is described as follows: it is the $A$-semilinear action on $\Mod_A(\mathcal C) \simeq \Mod_{M_k(A)}(\mathcal C)$ given by the $G$-equivariant map
    \[ BG \xrightarrow{\chi} B\Z/k \xrightarrow{\widehat{u}} B\PGL_k(A) = B\Aut_{A}(M_k(A)) = B\Aut_{A}(\Mod_A(\mathcal C)). \]
This is the same action that Baker, Richter and Szymik define at the level of $M_k(A)$.
\end{proof}



\begin{thm} \label{lem:Galois symbol for ring spectra}
Suppose given a strict unit $u \neq 1 \in \pi_0 \G_m(\bm 1)$. Its image in $\pi_0 \GL_1(\bm 1)$ is detected in the HFPSS for $\pic(\mathcal C) \simeq \pic(A)^{hG}$ by a class $v \in E_2^{s,s+1}$, and we assume that one of the following holds:
\begin{enumerate}
    \item $v$ is in positive filtration;
    \item $v$ is in filtration zero, and has nonzero image in $\widehat{H}^0(G, (\pi_0 A)^\times)$.
\end{enumerate}

Then the cyclic algebra $(A, \chi, u)$ is detected by the symbol
    \[ \beta(\chi) \cup v \in H^{s+2}(G, \pi_{s+1} \pic(A)) = E_2^{s+2,s+1}. \]
In particular, $\beta(\chi) \cup v$ is a permanent cycle. If it survives to $E_\infty$ then $(A, \chi, u) \neq \bm 1 \in \Br(\bm 1 \mid A)$.
\end{thm}

\begin{proof}
The square of $G$-spaces \cref{diag:square of spaces for cyclic algebra formula} gives rise to a commutative square of HFPSS as below:
\[\begin{tikzcd}[column sep=scriptsize, column sep = small]
	{H^*(G, \Z/k)} & {H^*(G,\Z)} & {\pi_* (B\Z/k)^{BG}} & {\pi_* (B^2\Z)^{BG}} \\
	{H^*(G, \pi_* \mathrm{PGL}_k A)} & {H^*(G, \pi_* \Picc(A))} & {\pi_* (B \mathrm{PGL}_k A)^{hG}} & {\pi_* (B\Picc(A))^{hG}}
	\arrow[""{name=0, anchor=center, inner sep=0}, "{(B\widehat{u})^{hG}}"', from=1-3, to=2-3]
	\arrow["{\beta^{BG}}", from=1-3, to=1-4]
	\arrow["{(B^2u)^{hG}}", from=1-4, to=2-4]
	\arrow["{\delta^{hG}}"', from=2-3, to=2-4]
	\arrow["{\beta_*}", from=1-1, to=1-2]
	\arrow[""{name=1, anchor=center, inner sep=0}, "{u_*}", from=1-2, to=2-2]
	\arrow["{\widehat{u}_*}"', from=1-1, to=2-1]
	\arrow["{\delta_*}"', from=2-1, to=2-2]
	\arrow[shorten <=40pt, shorten >=40pt, Rightarrow, from=1, to=0]
\end{tikzcd} \numberthis \label{diag:square of HFPSS for cyclic algebra formula} \]
Note that the maps on $E_2$ need not preserve filtration. By definition, $(A, \chi, u)$ is the image of $\chi \in \pi_0 (B\Z/k)^{BG}$ under the composite to $\pi_0 (B\Picc(A))^{hG}$, and is therefore detected on the $E_2$-page for $(B \Picc(A))^{hG}$ by $u_* \beta_* (\chi)$, as long as this class is nonzero. It is standard that the map $\beta_*$ is indeed the Bockstein, and we claim that the map
    \[ u_* : H^2(G, \Z) \to H^s(G, \pi_{s-1} \Picc(A)) \]
induced on $E_2$-pages by $B^2 u \colon B^2 \Z \to B \Picc(A)$ agrees with $v \cup -$. Indeed, this map can be identified with the composition
\begin{align*}
    \Map_G(EG_+, B^2 \Z) \simeq * \times \Map_G(EG_+, B^2 \Z) \xrightarrow{B^2 u \times id} & \Map_G(B^2\Z, B \Picc(A)) \times \Map_G(EG_+, B^2 \Z) \\
    \xrightarrow{\circ} & \Map_G(EG_+, B \Picc(A)),
\end{align*}
and the induced map on HFPSS shows that
    \[ u_* \beta(\chi) = [v \circ \beta(\chi)] = [ \beta(\chi) \cup v] \]
by compatibility of the cup and composition products. The class $\beta(\chi) \cup v$ is nonzero since Tate cohomology is $\beta(\chi)$-periodic, which gives the result.
\end{proof}

\subsection{Proof of \texorpdfstring{\cref{lem:widehat u is E1}}{Proposition~4.9}} \label{sec:lifting strict units}
We will write $\CMon \subset \Fun(\Fin_*, \mathcal S)$ for the \inftycat{} of special $\Gamma$-spaces, which admits a unique symmetric monoidal structure making the free commutative monoid functor $\mathcal S \to \CMon$ symmetric monoidal \cite{ggn}; its unit is the nerve of the category of finite pointed sets and isomorphisms,
    \[ \mathcal F \coloneqq |N \mathrm{Fin}_*^\sim| \simeq \coprod_{\Z[\geq 0]} B \Sigma_n. \]

\begin{defn}
We will denote by $\Rig \coloneqq \Alg_{\mathbb E_1}(\CMon)$ the \inftycat{} of \emph{associative semirings}\footnote{that is, rings without \emph{n}egatives}. Likewise, we will write $\CRig \coloneqq \Alg_{\mathbb E_\infty}(\CMon)$ for the \inftycat{} of \emph{commutative semirings}. Since the adjunction $\mathcal S \rightleftarrows \CMon$ is symmetric monoidal, it passes to algebras to yield monoid semiring functors,
    \[ \mathcal F[-]: \Alg_{\mathbb E_1}(\mathcal S) \rightleftarrows \Rig \qquad \text{and} \qquad \mathcal F[-]: \Alg_{\mathbb E_\infty}(\mathcal S) \rightleftarrows \CRig, \]
which factor the respective monoid algebra functors to spectra.
\end{defn}

\begin{remark}
The key observation for proving \cref{lem:widehat u is E1}, pointed out to us by Maxime Ramzi, is that for $R \in \CRig$ and $u \in \G_m(R)$, the map $\widehat u \in \pi_0 \Map(\Z/k, \PGL_k(R\grp))$ exists in $\pi_0 \Map(\Z/k, \PGL_k(R))$ before group completion. This is because the matrix $\widetilde{u}$ in \cref{eqn:widetilde u in GLk R} is defined without subtraction, as are all its powers. This allows us to make use of the fact that the representing object $\mathcal F[u^{\pm 1}] \coloneqq \mathcal F[\Z]$ for $\G_m$ is substantially simpler than $\mathbb S[u^{\pm 1}]$, as the following remark shows. 
\end{remark}

\begin{remark} \label{rem:F[Z]}
Consider the square
\[\begin{tikzcd}
	{\mathcal S} & \CMon \\
	{\Alg_{\mathcal O}} & {\Alg_{\mathcal O}(\CMon)}
	\arrow["L", shift left, from=1-1, to=1-2]
	\arrow["R", shift left, from=1-2, to=1-1]
	\arrow["{F'}"', shift right, from=1-2, to=2-2]
	\arrow["{U'}"', shift right, from=2-2, to=1-2]
	\arrow["{L' = \mathcal F[-]}", shift left, from=2-1, to=2-2]
	\arrow["{R'}", shift left, from=2-2, to=2-1]
	\arrow["F"', shift right, from=1-1, to=2-1]
	\arrow["U"', shift right, from=2-1, to=1-1]
\end{tikzcd}\]
of \cite[§7]{ggn}. By \cref{lem:free O algebra is an adjointable transformation} below, we have an equivalence
    \[ RU' \mathcal F[u^{\pm 1}] = RU'L' \mathbb Z \simeq RLU \Z \simeq RL \varinjlim_n [-n,n] \in \mathcal S.  \]
Since $\CMon$ is preadditive and $R \colon \CMon \to \mathcal S$ preserves sifted colimits, the underlying space of $\mathcal F[u^{\pm 1}]$ is $\varinjlim_n \prod_{[-n,n]} \mathcal F$. One therefore has a pullback square
\[\begin{tikzcd}[ampersand replacement=\&, column sep = small]
	{ \prod_i B\Sigma_{n_i}} \& {\mathcal F[u^{\pm 1}]} \\
	{\left\{ \sum_i n_i u^{t_i} \right\}} \& {\N[][u^{\pm 1}]}
	\arrow[from=2-1, to=2-2]
	\arrow[from=1-2, to=2-2, "{\pi_0}"]
	\arrow[from=1-1, to=2-1]
	\arrow[from=1-1, to=1-2]
\end{tikzcd} \numberthis \label{eqn:components of mathcal F[Z]} \]
for any basepoint of $\mathcal F[u^{\pm 1}]$.
\end{remark}

\begin{lem} \label{lem:free O algebra is an adjointable transformation}
Let $(L \dashv R) \colon \mathcal C \rightleftarrows \mathcal D$ be a symmetric monoidal adjunction between presentably symmetric monoidal \inftycats{}. Let $\lambda$ be an uncountable regular cardinal, and suppose that $\mathcal O^\otimes \to \mathrm{Sym}^\otimes$ is a fibration of $\infty$-operads compatible with colimits. Then the square
\[\begin{tikzcd}[column sep = small]
	{\mathcal C} & {\mathcal D} \\
	{\Alg_{\mathcal O}(\mathcal C)} & {\Alg_{\mathcal O}(\mathcal D)}
	\arrow["L", from=1-1, to=1-2]
	\arrow["{F'}", from=1-2, to=2-2]
	\arrow["{L'}"', from=2-1, to=2-2]
	\arrow["F", from=1-1, to=2-1]
\end{tikzcd}\]
is right adjointable.
\end{lem}

\begin{proof}
In the notation of \cref{rem:F[Z]}, we want to show that the mate
    \[ LU \to U'F'LU \simeq U'L'FU \to U'L' \label{eqn:mate for adjointability of AlgO} \numberthis \]
is an equivalence, and since $\Alg_{\mathcal O}(\mathcal C)$ is generated under sifted colimits by free algebras (as noted in the proof of \cite[Corollary~3.2.3.3]{ha}) and $U, U'$ preserve sifted colimits, it is enough to check this on free algebras $FC$, $C \in \mathcal C$. But $L'F \simeq F'L$, and using the description of $F$ as an operadic Kan extension \cite[Proposition~3.1.3.13]{ha} we factor \cref{eqn:mate for adjointability of AlgO} as an equivalence
    \[ LUFC = L \left( \coprod_{\N} \mathcal O(n) \otimes_{\Sigma_n} C^{\otimes n} \right) \xrightarrow{\sim} \coprod_{\N} \mathcal O(n) \otimes_{\Sigma_n} (LC)^{\otimes n} = U'F' LC \simeq U'L'FC. \qedhere \]
\end{proof}

Every nonzero coefficient in the matrix \cref{eqn:widetilde u in GLk R} is $1$, so it is a consequence of \cref{rem:F[Z]} that all relevant components of $\mathcal F[u^{\pm 1}]$ are contractible (which is certainly not true of $\mathbb S[\Z] = \mathbb S [u^{\pm 1}]$). This will allow us to prove that there is an essentially unique lift of semiring maps
\[\begin{tikzcd}[ampersand replacement=\&]
	{\Z/k} \& {\PGL_k(R)} \\
	\& {\pi_0 \PGL_k (R)}
	\arrow["{\widehat{u}}"', from=1-1, to=2-2]
	\arrow["{\widehat{u}}", dashed, from=1-1, to=1-2]
	\arrow[from=1-2, to=2-2]
\end{tikzcd}\]
which is functorial in $u$. In other words, we will deduce \cref{lem:widehat u is E1} from the following result:

\begin{prop} \label{lem:universal cyclic algebra is unique}
The fibre over $u \mapsto \widehat{u}$ of
    \[ \Map_{\mathcal P(\CRig\op)}(\G_m, \Map_{\mathbb E_1}(\Z/k, \PGL_k)) \to \Map_{\mathcal P(\CRig\op)}(\G_m, \Map_{\mathbb E_1}(\Z/k, \pi_0 \PGL_k)) \]
is contractible. 
\end{prop}

\begin{proof}
To prove the claim, note that $\G_m \in \mathcal P(\CRig\op)$ is corepresented by $\mathcal F[u^{\pm 1}]$, and so we are interested in the fibre of
    \[ \Map_{\mathbb E_1}(\Z/k, \PGL_k(\mathcal F[u^{\pm 1}])) \to \Map_{\mathbb E_1}(\Z/k, \pi_0 \PGL_k(\mathcal F[u^{\pm 1}])). \]
Using the inclusion $\Alg_{\mathbb E_1}(\mathcal S) \subset \Fun(\Delta \op, \mathcal S)$ \cite[Prop.~4.1.2.10]{ha} we have that
    \[ \Map_{\mathbb E_1}(\Z/k, \PGL_k(\mathcal F[u^{\pm 1}])) \simeq \lim_{[n] \in \Delta\op} \Map_{\mathcal S}(\Z/k^{\times n}, \PGL_k(\mathcal F[u^{\pm 1}])^{\times n}), \]
and it therefore suffices to prove that
each of the maps
    \[ \Map_{\mathcal S}(\Z/k^{\times n}, \PGL_k(\mathcal F[u^{\pm 1}])) \to \Map_{\mathcal S}(\Z/k^{\times n}, \pi_0 \PGL_k(\mathcal F[u^{\pm 1}])) \]
has contractible fibre over $(\widehat{u}^{j_1}, \dots, \widehat{u}^{j_k})$; for this, it is enough to show that
    \[ \PGL_k(\mathcal F[u^{\pm 1}]) \times_{\pi_0 \PGL_k(\mathcal F [u^{\pm 1}])} \{\widehat{u}^j\} \simeq * \]
for $j \geq 0$. But
\begin{align*}
    \PGL_k (\mathcal F[u^{\pm 1}]) &= \Aut_{\Alg_{\mathbb E_1}(\Mod_{\mathcal F[u^{\pm 1}]})}(M_k (\mathcal F[u^{\pm 1}])) \\
    & \subset \End_{\Alg_{\mathbb E_1}(\Mod_{\mathcal F[u^{\pm 1}]})}(M_k (\mathcal F[u^{\pm 1}])) \\
    &\simeq \lim_{n \in \Delta\op} \End_{\mathcal F[u^{\pm 1}]}(M_k (\mathcal F[u^{\pm 1}])^{\times n}),
\end{align*}
appealing this time to \cite[Prop.~4.1.2.11 and Theorem~2.3.3.23]{ha} for the final equivalence. Since $\Mod_{\mathcal F[u^{\pm 1}]}(\CMon)$ is semiadditive, we reduce to showing that the fibre
    \[ \End_{\mathcal F[u^{\pm 1}]}(M_k(\mathcal F[u^{\pm 1}])) \to \pi_0 \End_{\mathcal F[u^{\pm 1}]}(M_k(\mathcal F[u^{\pm 1}])) \]
over conjugation by $\widehat{u}^j$ is contractible. But any element in the right-hand side can be expressed as a sum of elementary matrix operations,
    \[ \sum_{1 \leq s,t,s',t' \leq k} a(s,t,s',t')(e_{st} \mapsto e_{s't'}), \]
where $a(s,t,s',t') = n(s,t,s',t') \, u^{r(s,t,s',t')} \in \pi_0 \mathcal F[u^{\pm 1}] = \mathbb N[u^{\pm 1}]$; the component of $\End_{\mathcal F[u^{\pm 1}]}(M_k (\mathcal F[u^{\pm 1}]))$ at this basepoint is
    \[ \prod_{s,t,s',t'} B \Sigma_{n(s,t,s',t')} \]
by \cref{eqn:components of mathcal F[Z]}. In particular, basic matrix algebra shows that the coefficients $n(s,t,s',t')$ for $\widehat{u}^j$ are all of the form $0$ or $1$, and so
    \[ \End_{\mathcal F[u^{\pm 1}]}(M_k(\mathcal F[u^{\pm 1}])) \times_{\pi_0 \End_{\mathcal F[u^{\pm 1}]}(M_k(\mathcal F[u^{\pm 1}]))} \{\widehat{u}^j\} \simeq \prod B\Sigma_0 \times \prod B\Sigma_1 \simeq *. \qedhere \]
\end{proof}
\begin{proof}[Proof (\cref{lem:widehat u is E1})]
We defined $\widehat u$ as an element of
    \[ \pi_0 \Map_{\mathcal P(\CRig \op)}(\G_m, \Map(\Z/k, \PGL_k)) = \pi_0 \Map(\Z/k, \PGL_k(\mathcal F[u^{\pm 1}])). \]
Write $P_k \coloneqq \PGL_k(\mathcal F[u^{\pm 1}])$. Forming the commutative square
\[\begin{tikzcd}[column sep=small]
	{\pi_0 \Map_{\mathbb E_1}(\Z/k,P_k)} & {\pi_0 \Map_{\mathbb E_1}(\Z/k,\pi_0P_k)} \\
	{\pi_0 \Map(\Z/k,P_k)} & {\pi_0 \Map(\Z/k,\pi_0 P_k)}
	\arrow[no head, equals, from=2-1, to=2-2]
	\arrow[from=1-2, to=2-2]
	\arrow[from=1-1, to=1-2]
	\arrow[from=1-1, to=2-1]
\end{tikzcd}\]
we observe that $\widehat{u}$ lifts to $\pi_0 \Map_{\mathbb E_1}(\Z/k, \pi_0 P_k) = \Hom_{\Ab}(\Z/k, \pi_0 P_k)$, by its definition. We thus obtain from \cref{lem:universal cyclic algebra is unique} a lift to $\pi_0 \Map_{\mathbb E_1}(\Z/k, P_k)$, and hence
    \[ (u \mapsto \widehat u) \in \Map_{\mathcal P(\CAlg\op)}(\G_m, \Map_{\mathbb E_1}(\Z/k, \PGL_k)). \numberthis \label{eqn:E1 cyclic algebra in Sp} \]
after passing to group completions. Now suppose $\mathcal C$ is any stable homotopy theory; then $\mathcal C$ is tensored over $\mathcal Sp$, and $\mathbb S[u^{\pm 1}] \otimes \bm 1 \simeq \bm 1[u^{\pm 1}]$. One therefore gets an $\mathbb E_1$ map
    \[ \PGL_k(\mathbb S[u^{\pm 1}]) = \Aut_{\Alg_{\mathbb E_1}}(M_k(\mathbb S[u^{\pm 1}])) \to \PGL_k(\bm 1[u^{\pm 1}]) = \Aut_{\Alg_{\mathbb E_1}(\mathcal C)}(M_k(\bm 1 [u^{\pm 1}])). \numberthis \label{eqn:E1 cyclic algebras in C} \]
Since $\bm 1[u^{\pm 1}] \in \CAlg(\mathcal C)$ corepresents $\G_m$, one obtains the desired map
    \[ (u \mapsto \widehat{u}) \in \Map_{\mathbb E_1}(\Z/k, \PGL_k(\bm 1[u^{\pm 1}])) = \Map_{\mathcal P(\CAlg(\mathcal C)\op)}(\G_m, \Map_{\mathbb E_1}(\Z/k, \PGL_k)) \]
as the image of \cref{eqn:E1 cyclic algebra in Sp} on applying $\Map_{\mathbb E_1}(\Z/k, -)$ to \cref{eqn:E1 cyclic algebras in C}.
\end{proof}

\section{Computing \texorpdfstring{$\Br_1^0$}{the relative Brauer group}}
\label{sec:azumaya}

We are now ready to complete the proofs of \cref{intro:Br1 even,intro:Br1 odd}. At odd primes, we will see that the cyclic algebra construction of \cref{sec:cyclic} gives all possible Brauer classes. On the other hand, when $p=2$ not all classes on the $E_\infty$-page of the descent spectral sequence will be detected in this way. In this case, we will first compute the relative Brauer group $\Br(KO_2 \mid KU_2)$.

\subsection{Odd primes}

When $p \geq 3$, \cref{fig:Picard spectral sequence p odd} shows that $\Br_1^0 \subset \mu_{p-1}$. In fact, we can deduce from \cref{lem:Galois symbol for ring spectra} that this inclusion is an equality:

\begin{thm} \label{lem:Brauer p odd}
Let $p$ be an odd prime, and choose $\chi \colon \mu_{p-1} \cong \Z/p-1$. There is an isomorphism
    \[ \mu_{p-1} \xrightarrow{\sim} \Br_1^0, \]
given by the cyclic algebra construction $\omega \mapsto (KU_p^{h(1 + p\Z[p])}, \chi, \omega)$.
\end{thm}

\begin{proof}
By \cref{lem:Picard sseq p odd} there is an inclusion $\Br_1^0 \subset \mu_{p-1}$. Moreover, there are no differentials out of the $0$-stem in \cref{fig:Picard spectral sequence p odd} by the computation of $\Pic_1$ \cite[§4]{profinitedescent}. Since $H^2(\mu_{p-1}, (\pi_0 B)^\times) \cong \mu_{p-1} \{ \beta(\chi) \cup \omega \}$, \cref{lem:Galois symbol for ring spectra} implies it is enough to show that the roots of unity in $(\pi_0 \bm 1_\K)^\times = \Z[p]^\times$ are strict. But we have a commuting square
\[\begin{tikzcd}[ampersand replacement=\&, row sep = large]
	{\pi_0 \G_m(\mathbb S_p)} \& {\pi_0 \G_m(\bm 1_\K)} \\
	{\pi_0 GL_1(\mathbb S_p)} \& {\pi_0 GL_1(\bm 1_\K)}
	\arrow[from=1-1, to=2-1]
	\arrow[no head, equals, from=2-1, to=2-2]
	\arrow[from=1-2, to=2-2]
	\arrow[from=1-1, to=1-2]
\end{tikzcd}\]
and at odd primes the roots of unity are strict in $\mathbb S_p$ by \cite[Theorem~A]{carmeli_strictunits}.
\end{proof}

\begin{remark}
In fact, \cref{lem:Brauer p odd} also follows from \cref{lem:Br = Br'}. Indeed, \cref{fig:Picard spectral sequence p odd} shows that
    \[ \Br'(\mathcal Sp_\K \mid KU_p) \cong H^2(\Z[p]^\times, \Z[p]^\times) \cong H^2(\mu_{p-1}, \Z[p]^\times). \]
In particular, this group is killed in the $\mu_{p-1}$-Galois extension $\bm 1_\K \to KU_p^{h(1 + p\Z[p])}$, since the group 
    \[ \Br'(KU_p^{h(1 + p\Z[p])} \mid KU_p) = \pi_0 B\Picc(KU_2)^{h(1 + p\Z[p])} \]
is concentrated in filtration $s \leq 1$ of the Picard spectral sequence for the $(1 + p\Z[p])$-action. Since the extension $\bm 1_\K \to KU_p^{h(1 + p\Z[p])}$ is finite, \cref{lem:Br = Br'} yields the second isomorphism below:
\begin{align*}
    \Br'(\mathcal Sp_\K \mid KU_p) &\cong \Br'(\mathcal Sp_\K \mid KU_p^{h(1 + p\Z[p])}) \\
    &\cong \Br(\mathcal Sp_\K \mid KU_p^{h(1 + p\Z[p])}) \subset \Br_1^0.
\end{align*}
\end{remark}


\subsection{Completed K-theory} \label{sec:K-theory}
In this section we use Galois descent to compute the Brauer group $\Br(KO_p \mid KU_p)$. This builds on the integral case computed by Gepner and Lawson \cite{gl}, and we will therefore also determine the completion maps 
    \[  \Br(KO \mid KU) \to \Br(KO_p \mid KU_p). \]
The computation for $p = 2$ will be important for our main computation: we will show that the relative Brauer classes of $KO_2$ descend to $\bm 1_\K$, which will help us determine the group $\Br_1^0$. Therefore, we start with the computation in this case:

\begin{thm} \label{lem:Brauer of KO2}
At the prime two we have
    \[ \Br(KO_2 \mid KU_2) \simeq \Z/4, \]
and the completion map from $\Br(KO \mid KU)$ is injective.
\end{thm}

Since the extension $KO_2 \to KU_2$ is finite, it follows by combining \cref{fig:Brauer spectral sequence for KO2} with \cref{lem:Br = Br'} that $|\Br(KO_2 \mid KU_2)| = 4$. To prove the theorem, we need to prove there is an extension, and we do this by reducing to computations of \'etale cohomology.

\begin{defn}
Recall the \emph{\'etale locally trivial Brauer group} $\LBr(KO_2) \subset \Br(KO_2)$ of \cite{antieau-meier-stojanoska}; more generally, Antieau, Meier and Stojanoska define
    \[ \LBr(R) \coloneqq \pi_0 \Gamma B\Picc_{\mathcal O_R} \]
for any commutative ring spectrum $R$, where $B\Picc_{\mathcal O_R}$ is the sheafification of $B\Picc(\mathcal O_R)$ on the \'etale site of $\Spet R$ \footnote{Note that this is equivalent to the \'etale site of $\Spet \pi_0 R$.}. Explicitly (\cite[Lemma~2.17]{antieau-meier-stojanoska}) this is the group of Brauer classes that are trivialised in some faithful \'etale extension $R \to R'$ in the sense of \cite[Definition~7.5.0.4]{sag}. Likewise, for an extension $R \to S$ of commutative ring spectra we write
    \[ \LBr(R \mid S) \coloneqq \ker (\LBr(R) \to \LBr(S)) \subset \Br(R \mid S). \]
\end{defn}

The group $\LBr(R)$ is sometimes more computationally tractable than $\Br(R)$: for example, one can often reduce to \'etale cohomology of $\Spet \pi_0 R$, which gives access to the standard cohomological toolkit. When $R = KO_p$, this allows us to use Gabber-Huber rigidity in the proof of \cref{lem:Brauer of KO2}.

\begin{remark}
In the setting of unlocalised \Einfty-rings, one always has compact generators and hence $\LBr(R) \cong \LBr'(R)$. One may (rightly) worry about the difference between the groups of unlocalised and of $K(1)$-local Brauer classes, since the results of \cite{antieau-meier-stojanoska} pertain to the former. While we do not know if the two groups agree in general (even for nice even-periodic rings), in our applications this is taken care of by restricting to \emph{relative} Brauer classes. Indeed, in that case we are computing the space
    \[ \Brr'(\bm 1 \mid A) = B\Picc(\Mod_A(\mathcal C))^{hG}, \]
and by \cite[Remark~3.7]{heard-mathew-stojanoska} the canonical map
    \[ \iota_A : \Picc(\Mod_A) \to \Picc(\Mod_A(\mathcal Sp_\K)) \]
is an equivalence of infinite loop-spaces in the following cases:
\begin{itemize}
    \item $A = E(k, \Gamma)$ is any Morava E-theory,
    
    \item $A$ admits a descendable extension $A \to B$ for which $\iota_B$ is an equivalence.
\end{itemize}
For example, this means that the two possible meanings of the expression $\LBr'(KO_2)$ agree. In fact, by \cref{lem:Br = Br'} we know that any element of $\LBr'(KO_2 \mid KU_2)$ lifts to an Azumaya algebra (in both the localised and unlocalised setting). This means that there is no ambiguity in writing $\LBr(KO_2 \mid KU_2)$ below.
\end{remark}

We begin with a preliminary computation; the following is essentially \cite[Proposition~3.8]{antieau-meier-stojanoska}:

\begin{prop} \label{lem:Picard sheaf on Spet Z2}
The \'etale sheaf $\pi_0 \pic(\mathcal O_{KO_2})$ fits in a nonsplit extension
    \[ 0 \to i_* \Z/2 \to \pi_0 \pic(\mathcal O_{KO_2}) \to \Z/4 \to 0, \]
where $i \colon \Spet \F[2] \to \Spet \Z[2]$ is the inclusion of the closed point. Moreover, $i^* \pi_0 \pic(\mathcal O_{KO_2}) \cong \Z/8$.
\end{prop}

\begin{proof}
We specify the necessary adjustments from the case of integral K-theory $KO$. Recall that Antieau, Meier and Stojanoska compute the sheaf $\pi_0 \pic(\mathcal O_{KO})$ on $\Spet KO = \Spet \Z$, using the sheaf-valued HFPSS for
    \[ \pic(\mathcal O_{KO}) \simeq \pic(\mathcal O_{KU})^{hC_2}, \]
which is \cite[Figure~1]{antieau-meier-stojanoska}. The same figure gives the HFPSS for
    \[ \pic(\mathcal O_{KO_2}) \simeq \pic(\mathcal O_{KU_2})^{hC_2}, \]
as long as one correctly interprets the symbols as in \cite[Table~1]{antieau-meier-stojanoska}, replacing $\mathcal O = \mathcal O_{\Z}$ with $\mathcal O_{\Z[2]}$. The proofs of \cite[Lemmas~3.5 and~3.6]{antieau-meier-stojanoska} go through verbatim to give the $0$-stem in the $E_\infty$-page, so that $\pi_0 \pic(\mathcal O_{KO_2})$ admits a filtration
\[\begin{tikzcd}[column sep=small]
	{F^2} & {F^1} & {\pi_0 \pic(\mathcal O_{KO_2})} \\
	{i_* \Z/2} & {\Z/2} & {\Z/2}
	\arrow[two heads, from=1-2, to=2-2]
	\arrow[hook, from=1-2, to=1-3]
	\arrow[two heads, from=1-3, to=2-3]
	\arrow[no head, equals, from=1-1, to=2-1]
	\arrow[hook, from=1-1, to=1-2]
\end{tikzcd}\]
Now the determination of the extensions follows as in \cite[Proposition~3.8]{antieau-meier-stojanoska}, by using the exact sequence
    \[ H^1(\Spet \Z[2], \G_m) = \Pic(\Z[2]) \to \Pic(KO_2) \to H^0(\Spet \Z[2], \pi_0 \pic(\mathcal O_{ KO_2})) \]
of \cite[Proposition~2.25]{antieau-meier-stojanoska}, and the fact that $\Z[2]$ is local so has trivial Picard group. In particular, since $i^*$ is exact we have that $i^* \pi_0 \pic(\mathcal O_{KO_2})$ admits a filtration by three copies of $\Z/2$, and a surjection from the constant sheaf $\Z/8$.
\end{proof}

\begin{proof}[Proof (\cref{lem:Brauer of KO2})]
By inspection of the HFPSS for the $C_2$ action on $\Picc(KU_2)$ (\cref{fig:Brauer spectral sequence for KO2}), the Brauer group is of order four. It remains to prove there is an extension relating these two classes. In fact, we will prove that
    \[ \Br(KO_2 \mid KU_2) \supset \LBr(KO_2 \mid KU_2) \cong \Z/4, \]
When $\pi_0 R$ is a regular complete local ring with finite residue field, the exact sequence \cite[Proposition~2.25]{antieau-meier-stojanoska} simplifies to an isomorphism
    \[ \LBr(R) \cong H^1(\Spet \pi_0 R, \pi_0 \pic(\mathcal O_R)), \]
since the cohomology of $\G_m$ vanishes \cite[1.7(a)]{mazur}. One has that $\pi_0 \pic(\mathcal O_{KU_2}) \simeq \Z/2$ is constant since $KU_2$ is even periodic with regular Noetherian $\pi_0$, while $\pi_0 \pic(\mathcal O_{KO_2})$ is torsion by \cref{lem:Picard sheaf on Spet Z2}. We can therefore use Gabber-Huber rigidity \cite{gabber,huber} to compute
\begin{align*}
    \LBr(KO_2) &\cong H^1(\Spet\Z[2], \pi_0 \pic(\mathcal O_{KO_2})) \cong H^1(\Spet \F[2], i^* \pi_0 \pic(\mathcal O_{KO_2})), \\
    \LBr(KU_2) &\cong H^1(\Spet\Z[2], \pi_0 \pic(\mathcal O_{KU_2})) \cong H^1(\Spet\F[2], i^* \pi_0 \pic(\mathcal O_{KU_2})).
\end{align*}
Since $i^* \pi_0 \pic(\mathcal O_{KO_2}) \simeq \Z/8$, we obtain
    \[ \LBr(KO_2) \cong \Z/8 \qquad \text{and} \qquad \LBr(KU_2) \cong \Z/2, \]
which implies that $\LBr(KO_2 \mid KU_2) \cong \Z/4$ or $\Z/8$; but $|\LBr(KO_2 \mid KU_2)| \leq 4$, so we are done.
\end{proof}

\begin{remark}
The generator of $\LBr(KO_2)$ is the Azumaya algebra $(KO_2\un, \Sigma KO_2\un)$ constructed in \cref{lem:KO2 algebra from H1 cocycle}, and the generator of $\LBr(KO_2 \mid KU_2)$ is $(KO_2\un, \Sigma^2 KO_2\un)$.
\end{remark}

\begin{remark}
In \cite[Lemma~V.3.1]{mqrt}, May proves a splitting of infinite loop-spaces
    \[ BU^{\otimes} \simeq \tau_{\leq 3} BU^{\otimes} \times \tau_{\geq 4} BU^{\otimes}, \]
where the monoidal structure is given by tensor product of vector spaces. The map $BU(1) = \tau_{\leq 3} BU \to BU$ is induced from a map at the level of bipermutative categories, and classifies the canonical line bundle. One can ask if this extends to a splitting of $\gl_1 KU$, and if this splitting happens $C_2$-equivariantly.

As an aside to \cref{lem:Brauer of KO2}, we deduce that the completed, equivariant analogue of this splitting fails:

\begin{cor}
There is no $C_2$-equivariant splitting
    \[ \gl_1 KU_2 \simeq \tau_{\leq 3} \gl_1 KU_2 \oplus \tau_{\geq 4} \gl_1 KU_2. \]
\end{cor}

\begin{proof}
The generators in the $(-1)$-stem of the $E_\infty$-page of \cref{fig:Brauer spectral sequence for KO2} are in filtrations two and six. If such a splitting did exist, there could be no extension between them in $\Br(KO_2 \mid KU_2)$.
\end{proof}

We do not know if $\gl_1 KU$ (or its $2$-completion) splits equivariantly: note that the obstruction in the case of $KU_2$ comes in the form of an extension on a class originating in $H^2(\Z[2]^\times, 1 + 4\Z[2]) \subset H^2(\Z[2]^\times, \Z[2]^\times)$.
\end{remark}

Let us briefly also mention the case when $p$ is odd; this will not be necessary for the computation of $\Br_1^0$ at odd primes.

\begin{prop}
When $p$ is odd, we have
    \[ \Br(KO_p \mid KU_p) \cong \Z/2, \]
and the map $\Br(KO \mid KU) \to \Br(KO_p \mid KU_p)$ is zero.
\end{prop}

\begin{figure}
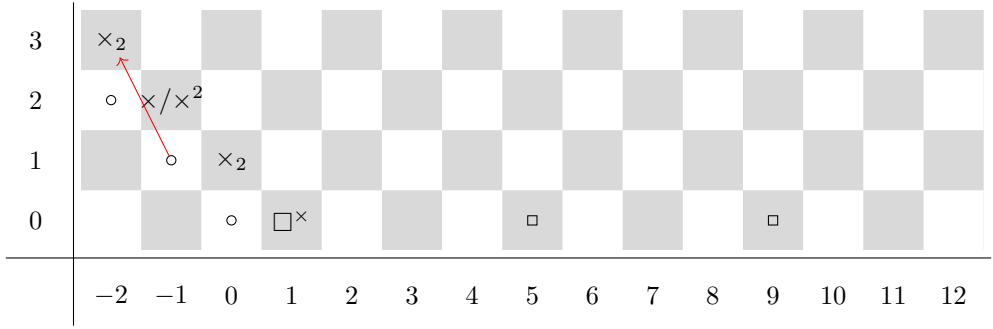

    \centering
    \printpage[grid = chess, xrange = {-2}{12}, yrange ={0}{3}, page = 2, name = BrauerKOp]
    \caption[HFPSS for $\pic(KO_p)$]{The HFPSS for $\pic(KO_p) \simeq \pic(KU_p)^{C_2}$ at odd primes. Here $\circ = \Z/2$, $\times = \mu_{p-1}$, $\times_2 = \mu_{p-1}[2] = C_2$ and $\square = \Z[p]$.}
    \label{fig:Picard spectral sequence for KO_p to KU_p p odd}
\end{figure}

\begin{proof}
Since $\Z[p]$ is local away from $2$, the $E_2$-page takes the form in \cref{fig:Picard spectral sequence for KO_p to KU_p p odd}, from which $\Br(KO_p \mid KU_p) \cong \mu_{p-1}/\mu_{p-1}^2 \cong \Z/2$ follows by \cref{lem:Br = Br'}. To describe the generators, note that the roots of unity $\mu_{p-1} \subset \pi_0 KO_p^\times$ are strict, since they are so in the $p$-complete sphere \cite{carmeli_strictunits}. Choosing $\chi: C_2 \cong \Z/2$, \cref{lem:Galois symbol for ring spectra} implies that the cyclic algebra construction
    \[ \omega \mapsto (KU_p, \chi, \omega) \]
yields $\beta(\chi) \cup \mu_{p-1} = H^2(C_{2}, \mu_{p-1}) \cong \Br(KO_2 \mid KU_2)$. The map from $\Br(KO \mid KU)$ is zero, since $\Br(KO_2 \mid KU_2)$ is detected in filtration $2$ and $\Br(KO \mid KU)$ in filtration $6$.
\end{proof}

\subsection{The case \texorpdfstring{$p = 2$}{p = 2}}
Putting together the work of the previous sections, we complete the computation of $\Br_1^0$ at the prime two.

\subsubsection{Descent from \texorpdfstring{$KO_2$}{KO2}.}

\Cref{lem:Brauer of KO2} yields
\begin{align*}
    \pi_t \Brr(KO_2 \mid KU_2) = \pi_t \Brr'(KO_2 \mid KU_2) = \left\{ \begin{array}{ll}
         \Z/4 \{(KO_2\un, \Sigma^2 KO_2\un)\} & t = 0 \\
         \Z/8 \{\Sigma KO_2 \} & t = 1 \\
         \Z[2]^\times & t = 2 \\
         \pi_{t-2} KO_2 & t \geq 3
    \end{array} \right.
\end{align*}
We will use this to compute the group $\Br(\bm 1_\K \mid KU_2)$ by Galois descent along the $\Z[2]$-Galois extension $\bm 1_\K \to KO_2$. Namely, we use the iterated fixed points formula
    \[ \Brr'(\bm 1_\K \mid KU_2) \simeq (B\Picc(KU_2)^{hC_2})^{h(1 + 4\Z[2])} = \Brr(KO_2 \mid KU_2)^{h(1 + 4\Z[2])}. \]
to form the descent spectral sequence
    \[ E_2^{s,t} = H^s(\Z[2], \pi_t \Brr(KO_2 \mid KU_2)) \implies \pi_{t-s} \Brr'(\bm 1_\K \mid KU_2). \]
Since $\Z[2]$ has cohomological dimension one, there is no room for differentials and the spectral sequence collapses immediately. To determine $\Br'(\mathcal Sp_\K \mid KU_2)$, what remains to compute is the following:
\begin{itemize}
    \item The group $E_2^{0,1} = \Br(KO_2 \mid KU_2)^{1 + 4\Z[2]}$,

    \item The extension between the groups $E_\infty^{0,-1} = E_2^{0,-1} \cong \Br(KO_2 \mid KU_2)^{1 + 4\Z[2]}$ and $E_\infty^{1,0} = E_2^{1,0} \cong \Z/8\{(KO_2, \Sigma KO_2) \}$.
\end{itemize}

This is achieved in the next couple of results, and the result is displayed in \cref{fig:sseq for 1midKO2}. Note that we have shifted degrees by one to match other figures, so that the relative Brauer group is still computed by the $(-1)$-stem.

\begin{figure}
    \centering
    \printpage[scale = 1.2, x range = {-2}{10}, y range = {0}{1}, grid = chess, page = \infty, name = Br(1midKO2)]
    \caption[Descent spectral sequence for $\Brr'(\bm 1_\K \mid KU_2) \simeq \Brr(KO_2 \mid KU_2)^{h(1 + 4\Z[2])}$]{The descent spectral sequence for $\Brr'(\bm 1_\K \mid KU_2) \simeq \Brr'(KO_2 \mid KU_2)^{h(1 + 4 \Z[2])}$. To match other figures, we have shifted everything in degree by one (so one may think of this as the spectral sequence for $\Sigma^{-1}\br'$). The extension in the $0$-stem is $4 \in \Ext(\Z/8, \Z[2]) \simeq \Z/8$, which gives $\pi_0 \Sigma^{-1} \br'(\bm 1_\K \mid KU_2) = \Pic_1 = \Z[2] \oplus \Z/4 \oplus \Z/2$.}
    \label{fig:sseq for 1midKO2}
\end{figure}

\begin{prop} \label{lem:Br10 to Brauer of KO2 is surjective}
We have $\Br(KO_2 \mid KU_2)^{1 + 4\Z[2]} = \Z/4$, so the map
    \[ \Br'(\bm 1_\K \mid KU_2) \to \Br(KO_2 \mid KU_2) \]
is surjective.
\end{prop}

\begin{proof}
It suffices to prove that $\psi^* (KO_2\un, \Sigma^2 KO_2\un) \simeq (KO_2\un, \Sigma^2 KO_2\un)$, where $\psi = \psi^\ell$ is the Adams operation for a topological generator $\ell \in 1 + 4\Z[2]$. By \cref{lem:Azumaya algebras from H1 classes}, this class is given by
    \[ \left( \varphi^{\Sigma KO_2\un}_! : \Mod_{KO_2\un} \to \varphi^* \Mod_{KO_2\un} \right) \in \operatorname{Eq} \left(id, \varphi^*: B\Picc(KO_2\un) \rightrightarrows B\Picc(KO_2\un) \right), \]
where $\varphi = KO_2 \otimes \varphi_2$ is the Frobenius on $KO_2\un = KO_2 \otimes_{\mathbb S} \mathbb{SW}$. In particular, note that $\psi \otimes \mathbb{SW}$ commutes with the $\varphi$. Thus the proposition follows from the square
\[\begin{tikzcd}[row sep=small]
	{\Mod_{KO_2\un}} & {\Mod_{KO_2\un}} & {\varphi^*\Mod_{KO_2\un}} \\
	&& {\psi^* \varphi^* \Mod_{KO_2\un}} \\
	{\psi^* \Mod_{KO_2\un}} & {\psi^* \Mod_{KO_2\un}} & {\varphi^* \psi^* \Mod_{KO_2\un}}
	\arrow["{\Sigma^2}", from=1-1, to=1-2]
	\arrow["{\varphi_!}", from=1-2, to=1-3]
	\arrow["{\psi_!}"', from=1-1, to=3-1]
	\arrow["{\Sigma^2}"', from=3-1, to=3-2]
	\arrow["{\varphi_!}"', from=3-2, to=3-3]
	\arrow["\sim", from=2-3, to=3-3]
	\arrow["{\psi_!}", from=1-3, to=2-3]
\end{tikzcd}\]
whose commutativity is witnessed by the natural equivalence
    \[ \psi_! \varphi_! \Sigma^2 \simeq \varphi_! \psi_! \Sigma^2 \simeq \varphi_! \Sigma^2 \psi_!. \qedhere \]
\end{proof}

\begin{prop} \label{lem:Br' p = 2}
The relative cohomological Brauer group at the prime two is
    \[ \Br'(\mathcal Sp_\K \mid KU_2) \cong \Z/8 \oplus \Z/4. \]
\end{prop}

\begin{proof}
Based on \cref{fig:sseq for 1midKO2}, what remains is to compute the extension from $\Br'(\bm 1_\K \mid KO_2) \cong \Z/8 \{(KO_2, \Sigma KO_2)\}$ to $\Br(KO_2 \mid KU_2)^{1  + 4\Z[2]} \cong \Z/4 \{(KO_2\un, \Sigma^2 KO_2\un)\}$. We will conclude by showing that the extension is split.

To prove the claim, note that both $(KO_2, \Sigma KO_2)$ and $(KO_2\un, \Sigma^2 KO_2\un)$ split over $KO_2\un$, so that the inclusion of $\Br'(\bm 1_\K \mid KU_2)$ in $\Br'(\bm 1_\K)$ factors as
\[\begin{tikzcd}[column sep = small]
	& {\Br'(\bm 1_\K \mid KO_2\un)} \\
	{\Br'(\bm 1_\K \mid KU_2)} & {\Br'(\bm 1_\K)}
	\arrow[hook, from=2-1, to=2-2]
	\arrow[hook, from=1-2, to=2-2]
	\arrow[dashed, from=2-1, to=1-2]
\end{tikzcd}\]
We will compute the relative Brauer group $\Br'(\bm 1_\K \mid KO_2\un)$ by means of the descent spectral sequence
    \[ H^s(\Z[2] \times \widehat{\Z}, \pi_t B\Picc(KO_2\un)) \implies \pi_{t-s} \Brr'(\bm 1_\K \mid KO_2\un), \]
which collapses at the $E_3$ page since $\Z[2] \times \widehat{\Z}$ has cohomological dimension two for profinite modules. To compute the $E_2$-page, note that
    \[ \Pic(KO_2\un) = \Z/8 \{ \Sigma KO_2\un \} \quad \text{and} \quad \pi_0 \GL_1(KO_2\un) = \W[]^\times, \]
by $C_2$-Galois descent from $KU_2\un \coloneqq E(\overline{\F}_2, \widehat{\G}_m) \simeq KU_2 \otimes_{\mathbb S} \mathbb{SW}$. The action on $\Pic(KO_2\un)$ is trivial, while the action on $\pi_0 \GL_1(KO_2\un)$ is trivial for the $\Z[2]$-factor, and Frobenius for the $\widehat{\Z}$-factor. In particular, note that $H^0(\widehat{\Z}, \W^\times) = \Z[2]^\times$; using the multiplicative splitting $\W^\times \simeq \overline{\F}_2^\times \times U_2$ (where $U_2 = \{x : \nu(x - 1) \geq 2\} \cong \W$), we see that
\begin{align*}
    H^1(\widehat{\Z}, \W^\times) &\cong H^1(\widehat{\Z}, \overline{\F}_2^\times) \oplus H^1(\widehat{\Z}, \W) \\
    & \cong 0 \oplus \varinjlim H^1(\Z/n, \W[](\F[2^n])) \\
    & = 0
\end{align*}
by Hilbert 90. This implies that $H^2(\Z[2] \times \widehat{\Z}, \pi_0 \GL_1(KO_2\un)) = 0$, since $\cd(\Z[2]) = \cd(\widehat{\Z}) = 1$ for profinite coefficients. The $(-1)$-stem of the $E_2$-page is hence concentrated in filtration one, and hence agrees with the $E_\infty$-page in this range. Thus
    \[ \Br'(\bm 1_\K \mid KO_2\un) = H^1(\Z[2] \times \widehat{\Z}, \Z/8) = \Z/8 \oplus \Z/8 \]
and $\Br'(\bm 1_\K \mid KU_2) \hookrightarrow \Z/8 \oplus \Z/8$, which implies the claim.
\end{proof}

\subsubsection{Generators at the prime two.}
To deduce \cref{intro:Br1 even} from \cref{lem:Br' p = 2} we will appeal to the results of \cref{sec:generators}.

\begin{thm} \label{lem:Brauer p = 2}
The relative Brauer group at the prime two is
    \[ \Br_0^1 = \Z/8 \oplus \Z/4. \]
\end{thm}

\begin{proof}
To lift the cohomological Brauer classes generating $\Br'(\mathcal Sp_\K \mid KU_2) \cong \Z/8 \oplus \Z/4$ to Azumaya algebras, it is enough by \cref{lem:Br = Br'} to prove that they are trivialised in some finite extension of $\bm 1_\K$. Recall from the previous subsection that:
\begin{enumerate}
    \item The generator of the $\Z/4$-factor is $(KO_2, \Sigma KO_2) \in \Br'(\mathcal Sp_\K \mid KO_2)$ (\cref{lem:KO2-trivial algebra from H1 cocycle}), and detected by $2 \in \Z/8 \cong H^1(\Z[2], \Pic(KO_2))$. Since $2$ is in the kernel of
        \[ H^1(\Z[2], \Pic(KO_2)) \to H^1(4\Z[2], \Pic(KO_2)), \]
    this cohomological Brauer class is trivialised in the $\Z/4$-Galois extension $KO_2^{h(1 + 16\Z[2])}$, so
        \[ (KO_2, \Sigma^2 KO_2) \in \Br'(\bm 1_\K \mid KO_2^{h(1 + 16\Z[2])}) \cong \Br(\bm 1_\K \mid KO_2^{h(1 + 16\Z[2])}) \subset \Br_1^0. \]

    \item Similarly, the generator of the $\Z/8$-factor is detected in the descent spectral sequence for $\bm 1_\K \to KO_2\un$ by $(1,0) \in \Z/8 \oplus \Z/8 \cong H^1(\widehat{\Z} \times \Z[2], \Z/8)$. We claim that this generator is trivialised in the extension $(KO_2\un)^{h(8\widehat{\Z} \times \Z[2])} \simeq \bm 1_\K \otimes_{\mathbb S} \mathbb{SW}_8$, where $\mathbb{SW}_8 \coloneqq W^+(\F[2^8])$. Indeed, since
        \[ (1,0) \in \ker \left( H^1(\widehat{\Z} \times \Z[2], \Z/8) \to H^1(8\widehat{\Z} \times \Z[2], \Z/8) \right) \]
    and the relative cohomological Brauer group $\Br'(\bm 1_K \otimes_{\mathbb S} \mathbb{SW}_8 \mid KO_2\un)$ is concentrated in filtration $s \leq 1$ of the Picard spectral sequence for $\bm 1_\K \otimes_{\mathbb S} \mathbb{SW}_8 \to KO_2\un$, we see that
        \[ [(1,0)] \in \Br'(\bm 1_\K \mid \bm 1_\K \otimes_{\mathbb S} \mathbb{SW}_8) \cong \Br(\bm 1_\K \mid \bm 1_\K \otimes_{\mathbb S} \mathbb{SW}_8) \subset \Br_1^0. \qedhere \]
\end{enumerate}
\end{proof}

\begin{remark}
Using \cref{lem:Br' p = 2}, we can also completely determine the behaviour of the Picard spectral sequence \cref{fig:Brauer spectral sequence p=2}. Recall the $E_2$-generators specified in \cref{rem:E2 generators}.
\begin{itemize}
    \item Under the base-change to $KO_2$, the generators $q_2$ and $q_6$ map to the $E_2$-classes representing $P_2, P_6 = P_2^2 \in \Br(KO_2 \mid KU_2)$. The splitting in \cref{lem:Brauer p = 2} of the surjection  $\Br_1^0 \twoheadrightarrow \Br(KO_2 \mid KU_2)$ of \cref{lem:Br10 to Brauer of KO2 is surjective} gives a canonical choice of classes $Q_2, Q_6 = Q_2^2 \in \Br_1^0$ lifting these. In particular, $q_2$ must also be a permanent cycle.

    \item Since $\Br_1^0 \supset \Br(\bm 1_\K \mid KO_2) \cong \Z/8$, the classes $q_1, q_2'$ and $q_4$ must also survive, and detect Brauer classes $Q_1, Q_2'$ and $Q_4$ trivialised over $KO_2$. We have $Q_2' = Q_1^2$ and $Q_4 = Q_1^4$ for this choice.
\end{itemize}

\begin{figure}[t]
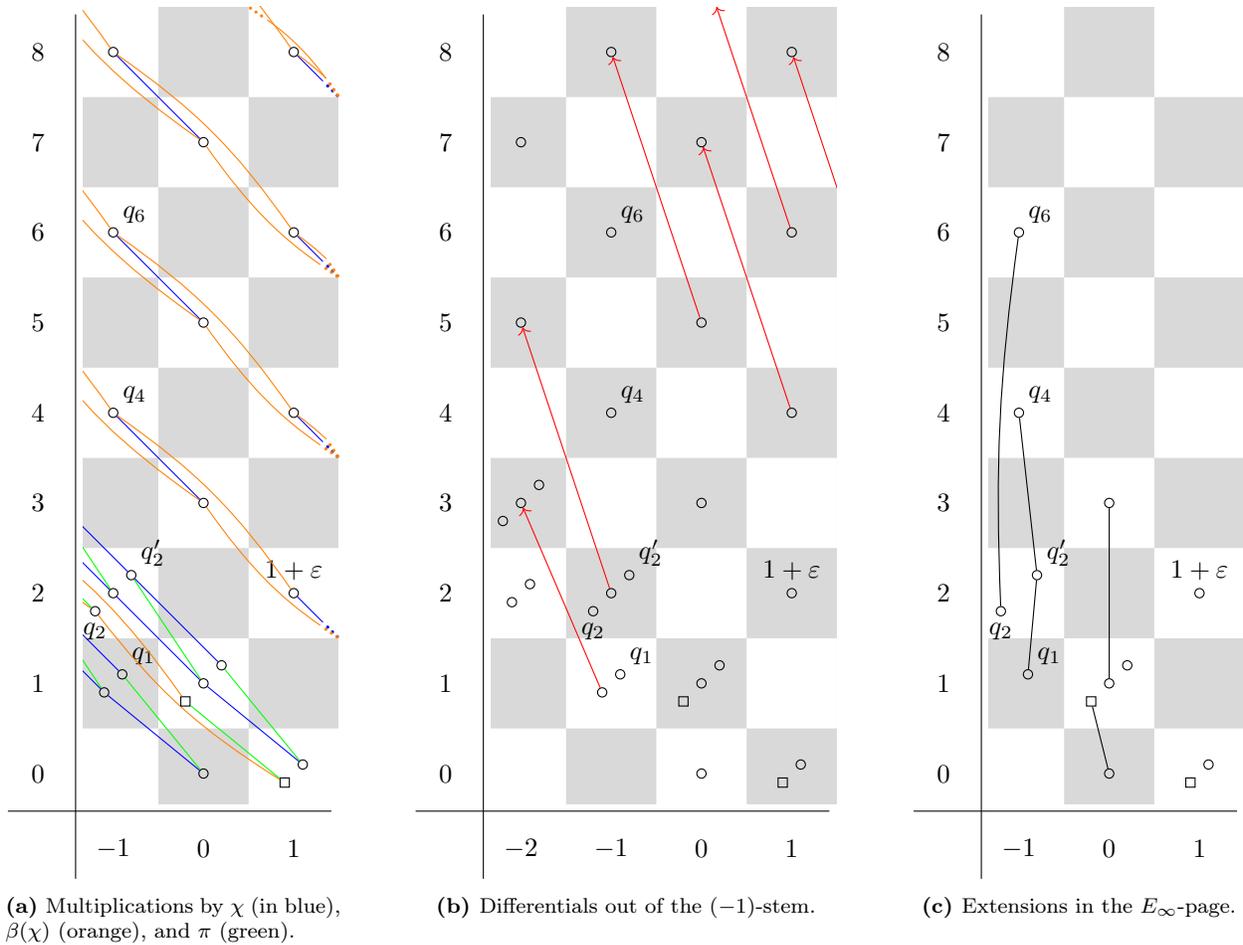

    \centering
    \begin{subfigure}{.27\textwidth}
    \centering
    \printpage[scale = 1.5, x range = {-1}{1}, y range = {0}{8}, grid = chess, page = 1, name = Brauersseven, draw struct lines, classes = {show name = none}]
    \caption{Multiplications by $\chi$ (in blue), $\beta(\chi)$ (orange), and $\pi$ (green).}
    \label{fig:module structure}
    \end{subfigure} \qquad %
    \begin{subfigure}{.36\textwidth}
    \centering
    \printpage[scale = 1.5, x range = {-2}{1}, y range = {0}{8}, grid = chess, page = 2--5, name = Brauersseven, classes = {show name = none}]
    \caption{Differentials out of the $(-1)$-stem. \\
    \phantom{a}}
    \end{subfigure} \qquad %
    \begin{subfigure}{.27\textwidth}
    \centering
    \printpage[scale = 1.5, x range = {-1}{1}, y range = {0}{8}, grid = chess, page = \infty, name = Brauersseven, classes = {show name = none}]
    \caption{Extensions in the $E_\infty$-page.\\
    \phantom{a}}
    \end{subfigure}
    \caption[Picard spectral sequence in low degrees]{Detailed view of the Picard spectral sequence (\cref{fig:Brauer spectral sequence p=2}) in low degrees.}
    \label{fig:Picard sseq detailed}
\end{figure}

In \cref{fig:Picard sseq detailed} we have displayed the $E_\infty$-page of the descent spectral sequence, including extensions. For the purposes of constructing explicit generators, we have also included the module structure over $H^*(\Z[2]^\times, \Z/2)$ and $H^*(C_2, \Z)$, as appropriate; in particular, in \cref{fig:module structure} we display multiplications by the generators
\begin{itemize}
    \item $\chi \in H^1(C_2, \Z/2) \cong \Hom(C_2, \Z/2)$,
    \item $\pi \in H^1(1 + 4\Z[2], \Z/2) \cong \Hom(\Z[2], \Z/2)$,
    \item $\beta(\chi) \in H^2(C_2, \Z)$.
\end{itemize}
\end{remark}

\begin{remark}
From the form of the spectral sequence, it follows that the class in $E_2^{7,5}$ survives to $E_\infty$---this should have implications for the nonconnective Brauer spectrum of $\mathcal Sp_\K$, as defined in \cite{haugseng_morita}.    
\end{remark}

Finally, we consider the consequences of \cref{lem:Galois symbol for ring spectra} at the prime two. The units $\Z/4^\times \subset \Z[2]^\times \subset \pi_0 \bm 1_\K^\times$ are not strict: for example, they are not strict in Morava E-theory by \cite[Theorem~8.17]{nullstellensatz}. In fact, we expect that the descent spectral sequence for $\mathbb G_m( \bm 1_\K) \simeq \mathbb G_m(E(\bar{\F}_2, \widehat{\G}_m))^{h (\widehat{\Z} \times \Z[2]^\times)}$ will yield
    \[ \pi_0 \mathbb G_m(\bm 1_\K) \cong \Z/2 \left\{ 1 + \varepsilon \right\} \subset (\pi_0 \bm 1_\K)^\times = (\Z[2][\varepsilon]/(2 \varepsilon, \varepsilon^2))^\times. \]
This will be discussed in future work. Nevertheless, we have the following corollary of \cref{lem:Galois symbol for ring spectra}:

\begin{cor}
For any $\chi \colon C_2 \cong \Z/2$, we have
    \[ Q_1^4 = Q_4 \coloneqq (KU^{h(1 + 4\Z[2])}, \chi, 1 + \varepsilon) \in \Br_1^0. \]
\end{cor}

\begin{proof}
The unit $1 + \varepsilon$ is strict by \cite[Corollary~5.5.5]{carmeli-yuan}, so the result follows from \cref{lem:Galois symbol for ring spectra} since $q_4 = \beta(\chi) \cup (1 + \varepsilon)$.
\end{proof}

\addcontentsline{toc}{section}{References}
\printbibliography

\end{document}